\documentclass[11pt]{amsart}
\usepackage{amsmath, amssymb, amsthm, hyperref, graphicx, scalerel}
\usepackage[margin=1.5in, top=1.4in, bottom=1.4in]{geometry}

\newtheorem{thm}{Theorem} \newtheorem*{thm*}{Theorem}
\newtheorem{prop}[thm]{Proposition} \newtheorem*{prop*}{Proposition}
\newtheorem{quest}[thm]{Question} \newtheorem*{quest*}{Question}
\newtheorem{cor}[thm]{Corollary} \newtheorem*{cor*}{Corollary}
\newtheorem{lemma}[thm]{Lemma} \newtheorem*{lemma*}{Lemma}
 \newtheorem*{conj*}{Conjecture}

 \newtheorem*{dfn*}{Definition}
\theoremstyle{remark}  \newtheorem*{example*}{Example}
  \newtheorem*{rmk*}{Remark}
\newtheorem*{ack*}{Acknowledgements}

\begin{document}

\title{Patterns in Random Permutations}

\pagestyle{plain}

\author{Chaim Even-Zohar}

\address{Department of Mathematics, 
UC Davis, 
One Shields Ave, 
Davis CA 95616}
\email{ chaim@math.ucdavis.edu }

\address{
The Alan Turing Institute, London, NW1 2DB, UK}
\email{ cevenzohar@turing.ac.uk }

\date{}

\begin{abstract}
Every $k$ entries in a permutation can have one of $k!$ different relative orders, called patterns.  How many times does each pattern occur in a large random permutation of size $n$? 

The distribution of this $k!$-dimensional vector of pattern densities was studied by Janson, Nakamura, and Zeilberger (2015). Their analysis showed that some component of this vector is asymptotically multi-normal of order $1/\sqrt{n}$, while the orthogonal component is smaller. 

Using representations of the symmetric group, and the theory of U-statistics, we refine the analysis of this distribution. We show that it decomposes into $k$ asymptotically uncorrelated components of different orders in $n$, that correspond to $S_k$-representations. 

Some combinations of pattern densities that arise in this decomposition have interpretations as practical nonparametric statistical tests.
\end{abstract}


%\keywords{????, ????}

\maketitle

\section{Introduction}
\label{intro}

Consider a given permutation $\pi \in S_n$ and a fixed $k \leq n$. We look at all the $\tbinom{n}{k}$ restrictions of $\pi$ to $k$ entries, $\pi(a_1),\dots,\pi(a_k)$ where $a_1 < a_2 < \dots < a_k$. The relative ordering of such $k$ values induces a \emph{pattern} $\sigma \in S_k$ as follows. The $k$-pattern $\sigma$ is the unique permutation in $S_k$ such that $\sigma(i) < \sigma(j)$ iff $\pi(a_i) < \pi(a_j)$ for every $i$ and~$j$. 

\begin{example*}
The restriction of $\pi = \underline{4}\,\underline{1}\,2\,\underline{5}\,3$ to the marked entries induces the $3$-pattern~$\sigma = 2\,1\,3$.
\end{example*}

For each $\sigma \in S_k$ we denote by $N_{\sigma}(\pi)$ the number of times it occurs as a $k$-pattern in $\pi$. The \emph{density} of $\sigma$ in $\pi$ is defined as the proportion $P_{\sigma}(\pi) = N_{\sigma}(\pi) / \tbinom{n}{k}$, so that $P_{\sigma}(\pi) \in [0,1]$. The \emph{$k$-profile} of~$\pi$ is the $k!$-dimensional vector of all $k$-pattern densities, 
$$ \mathbf{P}_k(\pi) \;=\; (\,P_{\sigma}(\pi)\,)_{\sigma \in S_k} \;\in\; \mathbb{R}^{k!} $$

\begin{example*}
The 3-profile of $\pi = 4\,1\,2\,5\,3$ is $\mathbf{P}_3(\pi) = \left( \frac{2}{10}, \frac{2}{10}, \frac{2}{10}, \frac{1}{10}, \frac{3}{10}, \frac{0}{10} \right)$. 
\end{example*}

\begin{rmk*}
Here and in the following, we arbitrarily order the $k!$ coordinates by the \emph{lexicographic} order, e.g., $123, 132, 213, 231, 312, 321$, and $\mathbb{R}^{k!}$ is the vector space with these coordinates.
\end{rmk*}

Pattern densities in permutations give rise to extremal questions~\cite{price1997packing, albert2002packing, hasto2002packing, presutti2010packing, burstein2010packing, balogh2015minimum, sliacan2017improving}, and $k$-profiles play a role in the construction of limiting objects for permutations~\cite{kral2013quasirandom, hoppen2013limits, glebov2015finitely, kenyon2015permutations}, and in the context of permutation property testing \cite{hoppen2011testing, klimovsova2014hereditary}, and quasirandom permutations~\cite{cooper2004quasirandom, cooper2006permutation, kral2013quasirandom}. The case where some pattern densities vanish is studied extensively~\cite{bona2012combinatorics,  kitaev2011patterns, marcus2004excluded}. 

\medskip

This paper explores the pattern densities and the profiles of a \emph{typical} permutation. Let $\pi \in S_n$ be uniformly distributed, so that $\mathrm{Pr}[\pi = \pi_0] = \tfrac{1}{n!}$ for every $\pi_0 \in S_n$. We denote its densities and $k$-profile by the random variables $P_{\sigma n}$ and $\mathbf{P}_{kn}$ respectively, sometimes abbreviated to $P_{\sigma}$ and $\mathbf{P}_{k}$. 

For every $\sigma \in S_k$ the expected density satisfies $E[P_{\sigma n}] = \frac{1}{k!}$ regardless of~$n$. This readily follows from linearity of expectation when we sum over the $\tbinom{n}{k}$ restrictions of the permutation to $k$ entries. The variance satisfies $V[P_{\sigma n}] = O(1/n)$, shown by summing over two $k$-wise restrictions, and considering different cases of their intersection. This implies a law of large numbers in $\mathbb{R}^{k!}$ for the random profile. Namely, for every $k \in \mathbb{N}$, 
$$ \mathbf{P}_{kn} \; \xrightarrow{\;\: n \to \infty \:\;} \; \mathbf{u}_k \;:=\; \left(\tfrac{1}{k!},\dots,\tfrac{1}{k!}\right) $$
in probability.

What is the typical deviation of the pattern densities from this limit? This question was studied in a series of works~\cite{fulman2004stein, bona2007copies, bona2010three, burstein2010packing, janson2015asymptotic, hofer2017central}. It was shown that in fact $V[P_{\sigma n}] = \Theta(1/n)$ for every $\sigma \in S_k$. Furthermore, the normalized density satisfies a central limit theorem, 
$$ \sqrt{n} \left( P_{\sigma n} - \tfrac{1}{k!}\right) \; \xrightarrow{\;\: n \to \infty \:\;} \; N\left(0,C_{\sigma\sigma} \right) $$
in distribution, for some $C_{\sigma\sigma} > 0$.

Burstein and H\"{a}st\"{o}~\cite{burstein2010packing} studied the joint distribution of two pattern densities $P_{\sigma n}$ and $P_{\sigma' n}$ where $\sigma \in S_k$ and $\sigma' \in S_{k'}$. In general, two such densities can have a nonzero correlation, even asymptotically as $n \to \infty$. In their paper, they derived formulas for the limits
$$ C_{\sigma\sigma'} \;:=\; \lim_{n \to \infty} \mathrm{cov}\left[\sqrt{n} P_{\sigma n},\sqrt{n} P_{\sigma' n}\right] \;.$$
In vector notation, this yields the leading, order $1/n$ term of the covariance matrix of the random $k$-profile, 
$$ {C}_{k} \;:=\; \left( C_{\sigma\sigma'}\right)_{\sigma,\sigma' \in S_k} \;=\; \lim_{n \to \infty} n \:\mathrm{cov}\left[\mathbf{P}_{kn}\right] \;\in\; \mathbb{R}^{k! \times k!} \;.$$

Janson, Nakamura and Zeilberger~\cite{janson2015asymptotic}, 
established a multivariate central limit theorem for the normalized $k$-profile,
$$ \sqrt{n} \left( \mathbf{P}_{kn} - \mathbf{u}_{k}\right) \; \xrightarrow{\;\: n \to \infty \:\;} \; N\left(\mathbf{0},{C}_{k} \right) $$
in distribution in $\mathbb{R}^{k!}$. They also showed convergence of all joint moments. 

More details on the distribution of the profile and the structure of this limit are given in these two works. Most importantly, the limit distribution of the normalized $k$-profile is degenerate, being supported on a $(k-1)^2$-dimensional linear subspace of~$\mathbb{R}^{k!}$~\cite{janson2015asymptotic}. This means that the distance of the vector $\sqrt{n} \left( \mathbf{P}_{kn} - \mathbf{u}_{k}\right)$ from this subspace goes to zero in probability as $n \to \infty$. 

Consider all scalar projections of the profile, i.e.~linear combinations of pattern densities, $\sum_{ \sigma \in S_k} c_{\sigma} P_{\sigma n}$. Since the limit distribution is supported on a small subspace, the typical magnitude of many such functionals, those in the orthogonal complement, is not determined except for being~$o(1/\sqrt{n})$. 

This paper seeks to further investigate the profile's distribution, looking at all directions in the $k!$-dimensional space. We'll examine the different orders of magnitude in~$n$ that appear, and mention some applications where also the smaller-order statistics are significant.

\medskip

Our main result states that the distribution of the random $k$-profile fits well a certain decomposition of the linear space $\mathbb{R}^{k!}$ that comes from representation theory:
$$ \mathbb{R}^{k!} \;=\; V_0 \oplus V_1 \oplus \dots \oplus V_{k-1} $$
The subspace $V_r$ is the span of all matrix elements of the linear $S_k$ representations that correspond to partitions $\lambda \vdash k$ with $\lambda_1 = k-r$. The description of functions and distributions on $S_k$ according to representations is sometimes called their Fourier transform. See Section~\ref{back} for a detailed definition of~$V_r$. 

We use the standard inner product $\left\langle \mathbf{u},\mathbf{v} \right\rangle = \sum_{\sigma \in S_k}u_{\sigma}v_{\sigma}$ in $\mathbb{R}^{k!}$. It is known that the $k$ components $V_0,\dots,V_{k-1}$ are mutually orthogonal with respect to this inner product.

The following theorem shows that the typical order of magnitude of the profile in different directions in $\mathbb{R}^{k!}$ varies according to these $k$ subspaces. As usual, the notation $a_n \sim b_n$ stands for $a_n/b_n \to 1$ as $n \to \infty$.
\begin{thm}
\label{proj1}
Let $r < k$ and $\mathbf{v} \in V_r \setminus \{\mathbf{0}\}$. 
$$ E\left[\; \left\langle \mathbf{v}, \mathbf{P}_{kn} \right\rangle^2 \;\right] \;\sim\; \frac{C_{\mathbf{v}}}{n^r}  $$
for some $C_{\mathbf{v}}>0$.
\end{thm}

We note that for $r \neq 0$ the expression on the left hand side of Theorem~\ref{proj1} is in fact $V\left[\left\langle \mathbf{v}, \mathbf{P}_{kn} \right\rangle \right]$, the variance of that linear functional of the random $k$-profile. This follows from $E[\mathbf{P}_{kn}] = \mathbf{u}_k \in V_0$, so that $E\left[\left\langle \mathbf{v}, \mathbf{P}_{kn} \right\rangle\right]  = 0$ if~$\mathbf{v} \in V_1 \oplus \dots \oplus V_{k-1}$, by the orthogonality of these subspaces.

Denote by $\Pi_{r}:\mathbb{R}^{k!} \to V_r$ the orthogonal projection induced by the direct sum. We write the $k$-profile as
$$ \mathbf{P}_{kn} \;=\; \Pi_0 \mathbf{P}_{kn} + \Pi_1 \mathbf{P}_{kn} + \dots + \Pi_{k-1} \mathbf{P}_{kn} $$
and call $\Pi_r \mathbf{P}_{kn}$ the \emph{order $r$ component} of the $k$-profile. Theorem~\ref{proj1} determines its order of magnitude, 
$$ E \left[\; \left\| \Pi_r \mathbf{P}_{kn} \right\|^2 \;\right] \;\sim\; \frac{C'_r}{n^r} \;\;\;\;\;\;\;\;\;\; r \in \{0,1,\dots,k-1\} $$
for some~$C'_r>0$. The decomposition $V_0 \oplus V_1 \oplus \dots \oplus V_{k-1}$ is actually characterized by the theorem in the following sense. Given $V_0,\dots,V_{r-1}$ the subspace $V_r$ is the smallest possible among all decompositions that satisfy such $k$ asymptotic relations.

\medskip

Another feature of this decomposition is that the $k$ normalized components, $n^{r/2} \,\Pi_r \mathbf{P}_{kn}$, are \emph{asymptotically uncorrelated}. This is stated in the next theorem in terms of the cross-covariance matrix of two such vectors.

\begin{thm}\label{cross}
For every $r < s < k$
$$ E \left[\;\left( n^{r/2} \,\Pi_r \mathbf{P}_{kn} \right) \left( n^{s/2} \,\Pi_s \mathbf{P}_{kn} \right)^T \;\right] \; \xrightarrow{\: n \to \infty \: } \; \mathbf{0} $$
in the normed space $\mathbb{R}^{k! \times k!}$ .
\end{thm}

To demonstrate the above results we briefly describe some cases of small~$k$, and mention their applications to nonparametric statistics. See Section~\ref{apps} for more details on these special cases.

\begin{example*}
In the case $k=2$, the decomposition is $\mathbb{R}^{2!} = V_0 \oplus V_1$, where $V_0 = \mathrm{span}\{\tbinom{1}{1}\}$ and $V_1 = \mathrm{span}\{\tbinom{+1}{-1}\}$. 

Indeed $P_{12}+P_{21}=1$, a constant, while the order of $P_{12}-P_{21}$ is $n^{-1/2}$. The latter is a well known property of the statistical test \emph{Kendall's}~$\tau$ \cite{kendall1938new}, or equivalently the \emph{inversion number} of a random permutation.
\end{example*}

\begin{example*}
For $k=3$, the components of $\mathbb{R}^{3!} = V_0 \oplus V_1 \oplus V_2$ have dimensions $1,4,1$ respectively. 

As usual the deterministic relation $\sum_{\sigma} P_{\sigma n} = 1$ yields the constant component in $V_0$. The subspace $V_1$ corresponds to an asymptotically 4-normal component of order $n^{-1/2}$. One of the four principal axes of this distribution can be used to express Spearman's~$\rho$ correlation test~\cite{spearman1904proof}. Finally, $V_2$ yields the \emph{sign} projection 
$$ P_{123}-P_{132}-P_{213}+P_{231}+P_{312}-P_{321} $$ 
whose order of magnitude is only $1/n$. This lower order statistic was suggested as a circular correlation test by Fisher and Lee~\cite{fisher1982nonparametric}, who also derived its non-normal limit distribution. See also \cite{zeilberger2016doron} and \cite[Remark~4.7]{janson2015asymptotic}.
\end{example*}

\begin{example*} The decomposition for $k=4$ is $\mathbb{R}^{4!} = V_0 \oplus V_1 \oplus V_2 \oplus V_3$, of dimensions $1,9,13,1$, and orders of magnitude $1,\,n^{-0.5}, n^{-1}, n^{-1.5}$ respectively. 

Some order-$1/n$ projections coming from the subspace $V_2$ for $k \geq 4$ yield variants of Hoeffding's independence test for paired samples~\cite{hoeffding1948non, blum1961distribution, bergsma2014consistent}.
\end{example*}

The main tools in the proofs of Theorems~\ref{proj1} and \ref{cross} are the theory of \mbox{U-statistics}, and the construction of $S_k$-representations via Young symmetrizers. Section~\ref{back} contains some background on these topics. The proofs are presented in Section~\ref{proofs}.

Theorems~\ref{proj1}-\ref{cross} decompose the $k$-profile into $k$ orthogonal components, that are asymptotically pairwise-uncorrelated. One may proceed by decomposing every normalized component $n^{r/2}\Pi_r \mathbf{P}_{kn}$ according to an appropriate orthogonal basis of~$V_r$, so that the covariance matrix of all the resulting one-dimensional components is asymptotically diagonal. Such a procedure is sometimes called \emph{decorrelation} or \emph{principal component analysis} (PCA). In Section~\ref{diag} we describe such a diagonalization explicitly for every $k \leq 6$, which is again given by matrix elements of representations of~$S_k$. These preliminary results lead to interesting directions for future work on the distribution of the $k$-profile. 

Finally, in Section~\ref{apps} we discuss the relations and applications to nonparametric statistical tests, and quasirandomness. 

\section{Background}
\label{back}

We first review the representation theory of the symmetric group $S_k$, and define the decomposition of the profile in detail. Then we review the theory of U-statistics, to
 be applied in the analysis of the resulting components.

\bigskip \noindent
\textbf{Representations of the Symmetric Group}
\nopagebreak \medskip \\ \noindent
We present here some standard facts about group representations, and refer to Fulton and Harris~\cite[Lectures 1--4]{fulton1991representation} for a full exposition. We work over~$\mathbb{R}$ rather than~$\mathbb{C}$, since all the representations of $S_k$ over $\mathbb{C}$ are real.

A \emph{$d$-dimensional real representation} of a finite group $G$ is a map $\rho$ from $G$ to the linear group of $\mathbb{R}^d$, such that $\rho(g) \circ \rho(h) = \rho(gh)$ for every $g,h \in G$. If~there is no proper subspace $V$ of $\mathbb{R}^d$ such that $\rho(g)V=V$ for every $g \in G$, then the representation $\rho$ is called \emph{simple}, also known as \emph{irreducible}. Two representations $\rho,\rho'$ of $G$ are \emph{similar} if there exists a linear map $\tau$ such that $\rho'(g) = \tau^{-1} \circ \rho(g) \circ \tau$ for every~$g \in G$. Every finite group $G$ has a finite number of simple representations up to similarity.

The representation theory of the symmetric group is a classical well-studied area. There is a one-to-one correspondence between simple representations of $S_k$ and integer partitions of $k$. A~\emph{partition} $\lambda$ of $k \in \mathbb{N}$, denoted $\lambda \vdash k$, is a nonincreasing sequence of integers $\lambda_1 \geq \lambda_2 \geq \cdots \geq \lambda_{\ell} \geq 1$ such that $\lambda_1 + \dots + \lambda_{\ell} = k$. Let $\rho^{\lambda}$ denote the representation that corresponds to~$\lambda$, defined up to similarity. Also, let $d_{\lambda}$ denote the dimension of $\rho^{\lambda}$. These dimensions are known to satisfy $\sum_{\lambda \vdash k} d_{\lambda}^2 = k!$. See the proof of Lemma~\ref{cosets} in Section~\ref{proofs} for a brief construction of $\rho^\lambda$, or~\cite[\S 4.1]{fulton1991representation}.

Let $\lambda \vdash k$, and suppose that the simple representation~$\rho^\lambda$ of $S_k$ is given as $d_\lambda$-by-$d_\lambda$ matrices $R^\lambda(\sigma)$ for all $\sigma \in S_k$. For fixed $i,j \in \{1,\dots,d_{\lambda}\}$ consider the vector consisting of the $(i,j)$-entry of each matrix,
$$ \mathbf{R}^{\lambda}_{ij} \;=\; \left(\; R^{\lambda}_{ij}(\sigma) \;\right)_{\sigma \in S_k} \in\; \mathbb{R}^{k!} $$  
Such a vector is called a \emph{matrix element} of $\rho^{\lambda}$. The matrix elements of a simple real representation of $S_k$ span a $d_\lambda^2$-dimensional subspace of $\mathbb{R}^{k!}$. This subspace only depends on $\lambda$ since it is preserved under similarity. The matrix elements of all simple representations together span the whole space~$\mathbb{R}^{k!}$.

\bigskip \noindent
\textbf{The Decomposition $V_0 \oplus \dots \oplus V_{k-1}$}
\nopagebreak \medskip \\ \noindent
We define the subspaces $V_r$ of Theorems~\ref{proj1} and~\ref{cross}, based on the above notions, and then describe some special cases as examples. The parameter $k \in \mathbb{N}$ is kept fixed, and suppressed to simplify the notation.
\begin{dfn*}
For $r \in \{0,1,\dots,k-1\}$ let
$$ V_r \;:=\; \mathrm{span}\,\left\{ \left. \; \mathbf{R}^{\lambda}_{ij} \;\right|\; \lambda \vdash k \;\text{ with }\; \lambda_1 = k-r, \;\; 1 \leq i,j \leq d_\lambda \right\}$$
\end{dfn*}
In words, $V_r$ is spanned by the matrix elements of all representations $\rho^{\lambda}$ with $\lambda_1 = k-r$. Since the definition of $V_0,\dots,V_{k-1}$ uses each representation once, these spaces are orthogonal to each other, and $V_0 \oplus \dots \oplus V_{k-1} = \mathbb{R}^{k!}$. We demonstrate this definition through some basic examples~\cite[1.3]{fulton1991representation}.

\begin{example*}
There is a single partition $\lambda \vdash k$ with $\lambda_1=k$, the one-part partition. It corresponds to the one-dimensional \emph{trivial representation} of $S_k$, $R^{k} (\sigma) = [1]$ for every $ \sigma \in S_k$. Thus $V_0$ is the one-dimensional subspace $\mathrm{span}\,\{(1,\dots,1)\}$ in~$\mathbb{R}^{k!}$.
\end{example*}

\begin{example*} The only partition of $k$ whose largest part $\lambda_1=1$, is $(1,\dots,1)$ with $k$ parts. It corresponds to $R^{1 \cdots 1} (\sigma) = [\mathrm{sign}\,\sigma]$, the \emph{alternating representation} of $S_k$. This yields the last component, $V_{k-1}=\mathrm{span}\,\left\{(\mathrm{sign}\,\sigma)_{\sigma \in S_k}\right\}$, which is one-dimensional too.
\end{example*}

\begin{example*}
With $\lambda_1=k-1$, there is a two-part partition $(k-1,1)$, which corresponds to the $(k-1)$-dimensional \emph{standard representation} of $S_k$. The $(k-1)^2$-dimensional span of its matrix elements yields the subspace $V_1$.

Here is an explicit description of the standard representation.
Let $A(\sigma)$ be the $k \times k$ permutation matrix of $\sigma$, that is $A_{ij}(\sigma) = \delta_{\sigma(i)j}$, and let $U$ be a $k \times (k-1)$ matrix whose columns are any orthonormal basis of $(1,\dots,1)^{\perp}$ in $\mathbb{R}^k$. Then the standard representation can be given by the orthogonal matrices $R^{(k-1)1}(\sigma) = U^T A(\sigma)\, U$.
\end{example*}

We note that these three cases are in fact the only ones where $V_r$ is defined by a single simple representation. If $k=4$ for example, then we have two partitions with $\lambda_1 = 2$, namely $(2,2)$ and $(2,1,1)$, which correspond to representations of $S_4$ of dimensions $2$ and $3$ respectively. In this case, $V_2$ is a subspace of $\mathbb{R}^{4!}$ of dimension $2^2+3^2=13$.

\bigskip \noindent
\textbf{Two-Sided Cosets}
\nopagebreak \medskip \\ \noindent
We present an equivalent description of the above decomposition, which will come in useful in the proofs of Theorems~\ref{proj1} and~\ref{cross}.

Let $1 \leq l \leq k$, and consider the subgroup $S_l \leq S_k$, which we identify with all permutations that fix the first $k-l$ elements, that is, 
$$ S_l \;:=\; \left\{\left. \tau \in S_k \; \right| \; \tau(i)=i \text{ for } i \leq k-l \right\}$$ 
A \emph{two-sided coset} of 
$S_l$ in~$S_k$ is a subset of the form $\alpha S_l \beta = \left\{\left. \alpha \tau \beta \;\right|\; \tau \in S_l \right\}$
for some $\alpha, \beta \in S_k$. Equivalently, a two-sided $S_l$-coset contains $l!$ permutations, with $k-l$ fixed values in $k-l$ fixed positions, and the remaining $l$ positions and values are matched in all possible ways.

The following lemma essentially  characterizes the subspaces $V_0, \dots, V_{k-1}$ in terms of two-sided cosets.

\begin{lemma}
\label{cosets}
For $l \in \{1,2,\dots,k\}$ 
$$ \mathrm{span}\,\left\{  \; \mathbf{R}^{\lambda}_{ij} \;\left|\; \begin{array}{c} \lambda \vdash k \\ \lambda_1 < l \\ 1 \leq i,j \leq d_\lambda \end{array} \right. \right\}
\;=\;
\left\{ \mathbf{v} \in \mathbb{R}^{k!} \;\left|\; \forall \alpha,\beta \in S_k: \;\sum\limits_{\sigma \in \alpha S_l \beta} v_\sigma = 0 \right. \right\} $$
\end{lemma}
The subspace on the left hand side is the sum of $V_r$ for all~\mbox{$r > k-l$}. This property thus characterizes $V_0,\dots,V_{k-1}$ as the unique orthogonal decomposition that is compatible with two-sided cosets in this fashion.

\begin{example*}
Let $k=3$. For $l=3$, the condition on cosets of $S_l$ boils down to the $6$ entries of $\mathbf{v}$ adding up to $0$, which leaves the $5$-dimensional space $V_1+V_2$. In the case $l=2$, every $S_2$-coset yields two entries of $\mathbf{v}$ that differ by a transposition in $S_3$, and have opposite signs. This indeed defines the $1$-dimensional space $V_2$. The condition for $l=1$ trivially yields $\{\mathbf{0}\}$ on both sides.
\end{example*}

Lemma~\ref{cosets} demonstrates the relation between Young subgroups of $S_k$ and its linear representations, in the special form needed here. We provide a proof of this lemma in Section~\ref{proofs}, based on 
standard results from Fulton and Harris~\cite{fulton1991representation}.

\bigskip \noindent
\textbf{U-Statistics}
\nopagebreak \medskip \\ \noindent
The $k$-profile's analysis involves certain summations over the $\tbinom{n}{k}$ potential occurrences of various $k$-patterns. Since we study the profile of a uniformly distributed permutation in $S_n$, such $\tbinom{n}{k}$ summands are identically distributed random variables. However, they aren't necessarily independent, and thus cannot be handled with the classical limit theorems for iid random variables. Instead, we employ the theory of U-statistics, a broader class of random sums that are suitable for this setting.

Here is the definition of U-statistics, coming from the foundational work of Hoeffding~\cite{hoeffding1948class}. Let $X_1, X_2, \dots$ be iid random variables, taking values in a measurable space~$\mathfrak{X}$. Let the \emph{kernel} function $f:\mathfrak{X}^k \to \mathbb{R}$ be a \emph{symmetric} measurable function, namely $f(x_1,\dots,x_k) = f(x_{\sigma(1)},\dots,x_{\sigma(k)})$ for every $\sigma \in S_k$. The following sequence of random variables, for $n \geq k$, are \mbox{\emph{U-statistics}} of order~$k$:
$$ U_n \;=\; \frac{1}{\binom{n}{k}} \sum\limits_{i_1< \dots< i_k} f\left(X_{i_1},\dots,X_{i_k}\right) $$
The sum is over ${\tbinom{n}{k}}$ terms, all the $k$-element subsets of $\{1,\dots,n\}$. We present here some basic properties of this class of random variables, and refer to~\cite{lee1990u, janson1997gaussian, korolyuk2013theory} for more.   

Given the kernel $f : \mathfrak{X}^k \to \mathbb{R}$ as above, we define a sequence of functions $f_0, f_1, \dots, f_k$, such that $f_i : \mathfrak{X}^i \to \mathbb{R}$.  
$$ f_i(x_1,\dots,x_i)\;:=\; E\left[f\left(x_1,\dots,x_i,X_{i+1},\dots,X_k\right)\right] $$
where the expectation is taken with respect to $X_{i+1},\dots,X_k$, independent and identically distributed as above. Note that $f_k=f$, while
$$ f_0 \;=\; E[f_i] \;=\; E[U_n]$$
for every $i \leq k$ and $n \geq k$. Here we shorthand $E[f_i]=E[f_i(X_1,\dots,X_i)]$ where $X_1,\dots,X_i$ are iid as above, and similarly for the variance. We assume throughout $V[f] < \infty$. By the law of total variance,
$$ 0 \;=\; V[f_0] \;\leq\; V[f_1] \;\leq\; \dots \;\leq\; V[f_k] \;=\; V[f]$$
The \emph{rank} of the function $f$, or the \emph{rank} of $U_n$, is defined to be the smallest~$r$ for which $V[f_r] \neq 0$. This number is sometimes named the \emph{degree of degeneracy}, and has a great impact on the properties of $U_n$.

The variance of a U-statistic can be expressed using the variance of the functions $f_r,\dots,f_k$ as follows:
$$ V[U_n] \;=\; \frac{1}{\binom{n}{k}} \sum\limits_{i=1}^k \binom{k}{i} \binom{n-k}{k-i} V[f_i] $$
This formula follows by splitting up a double sum over two $k$-element subsets, according to the size of the intersection~\cite[page 12]{lee1990u}. It immediately follows that the rank of a U-statistic determines its typical order of magnitude:
\begin{cor}
\label{rank}
Let $U_n$ be U-statistic with a kernel $f$ of order $k$ and rank $r$. Then 
$$ V[U_n] \;=\; \frac{\binom{k}{r}^2 r! \,V[f_r]}{n^r} + O\left(\frac{1}{n^{r+1}}\right)$$
\end{cor} 
Finally, we remark that many other properties of the distribution of $U_n$ depend on its rank $r$. For a U-statistic of rank one, for example, $\sqrt{n}U_n$ converges in distribution to a normal random variable. If the rank is two or more, then $n^{r/2} U_n$ is not asymptotically normal, and can have a limit distribution of various shapes, with the details depending on the kernel~$f$. However, its tail behavior always satisfies $\mathrm{Pr}\left[|n^{r/2}U_n| > t\right] = \exp\left(-\Theta\left(t^{2/r}\right)\right)$ as $t \to \infty$. In this work we focus on second moments, and do not discuss the distributions in more detail.

\section{Proofs}
\label{proofs}

We first make some observations about the random profile, which are relevant to the subsequent proofs.  

There is a standard construction relating a permutation in $S_n$ with the order statistics of $n$ points in the plane. Let $\{x_1,\dots,x_n\}$ be a set of $n$ points in $\mathbb{R}^2$. We use throughout the notation $x_i = (y_i,z_i)$ for the two coordinates, and consider  \emph{generic} sets, where all the $y$-coordinates are distinct, as well as the $z$-coordinates. The $y$-coordinates $y_1,\dots,y_n$ are written increasingly as \emph{order statistics}: $y_{(1)} < y_{(2)} < \dots < y_{(n)}$, and similarly $z_1,\dots,z_n$ are reordered as $z_{(1)} < z_{(2)} < \dots < z_{(n)}$. There is a unique permutation \mbox{$\pi \in S_n$} such that the points $\{x_1,\dots,x_n\}$ are given by $\{(y_{(i)}, z_{(\sigma(i))}) \;|\; 1 \leq i \leq n\}$. We call $\pi$ the permutation \emph{induced} by $x_1,\dots,x_n$, and denote it by  $\mathrm{perm}(x_1,\dots,x_n)$. Note that $\mathrm{perm}$ is symmetric with respect to shuffling its $n$ inputs.

Consider $\mathrm{perm}(X_1,\dots,X_n)$ where $X_1,\dots,X_n$ are independent random points, each one sampled from a continuous distribution $X = (Y,Z)$ in~$\mathbb{R}^2$. With probability one, such a point set is generic and the induced permutation $\pi \in S_n$ is well-defined. Observe that if $Y$ and $Z$ are independent, then $\pi$ is uniformly distributed. Similarly, every $k$ variables $X_{i_1},\dots,X_{i_k}$ induce a $k$-pattern of $\pi$, uniform in $S_k$.

In the subsequent analysis of the random profile~$\mathbf{P}_{kn}$ we assume it originates from such a random permutation $\mathrm{perm}(X_1,\dots,X_n)$ where the $X_i$ are independent and follow the continuous uniform distribution in the unit square~$[0,1]^2$.

Let $\mathbf{v} \in \mathbb{R}^{k!}$. Now we show as in~\cite{ janson2015asymptotic} that the scalar projection of the $k$-profile $\left\langle \mathbf{v}, \mathbf{P}_{kn} \right\rangle$ is a U-statistic of order $k$.
\begin{align*}
\left\langle \mathbf{v}, \mathbf{P}_{kn} \right\rangle
\;&=\;
\sum\limits_{\sigma \in S_k} v_{\sigma} P_{\sigma}(\mathrm{perm}(X_1,\dots,X_n)) 
\\
\;&=\; 
\sum\limits_{\sigma \in S_k} v_{\sigma}\,\frac{1}{\tbinom{n}{k}} \sum\limits_{i_1 < \dots < i_k} \begin{cases} 1 &\mathrm{perm}\left(X_{i_1},\dots,X_{i_k}\right) = \sigma \\ 0 & \text{otherwise} \end{cases}
\\
\;&=\; \frac{1}{\tbinom{n}{k}} \sum\limits_{i_1 < \dots < i_k} v_{\mathrm{perm}\left(X_{i_1},\dots,X_{i_k}\right)}
\end{align*}
This is indeed an order-$k$ U-statistic with the kernel 
$$ f(x_1,\dots,x_k) \;=\; v_{\mathrm{perm}\left(x_1,\dots,x_k\right)} $$ 
Note that although the entries of the induced $\mathrm{perm}(X_1,\dots,X_n)$ are dependent as random variables, the underlying $X_1,\dots,X_n$ are mutually independent.

\bigskip \noindent
\textbf{Proof of Theorem~\ref{proj1}}
\nopagebreak \medskip \\ \noindent
Let $\mathbf{v} \in V_r \setminus \{\mathbf{0}\}$ for some $r<k$. Recalling the definition in Section~\ref{back}, this means
$$ \mathbf{v} \;\in\; \mathrm{span}\,\left\{ \left. \; \mathbf{R}^{\lambda}_{ij} \;\right|\; \lambda \vdash k \;\text{with}\; \lambda_1 = k-r, \;\; 1 \leq i,j \leq d_\lambda \right\}$$

In the case $r=0$, the vector $\mathbf{v}$ is a nonzero multiple of the trivial  representation $\mathbf{R}^{k}_{11}$, hence $\mathbf{v} = (a,\dots,a)$ where $a \neq 0$. Then, for all~$n \geq k$ 
$$ E\left[\left\langle \mathbf{v}, \mathbf{P}_{kn} \right\rangle^2 \right] \;=\; E\left[\left( a\sum\nolimits_{\sigma \in S_k} P_{\sigma n} \right)^2 \right] \;=\; E\left[a^2\right] \;=\; a^2 $$
In particular, the order of this expectation is constant with respect to $n$, as claimed in the theorem. 

We hence pursue the proof assuming $0<r<k$. Since $\mathbf{v} \in V_r$, which is orthogonal to $\mathbf{1} \in V_0$, 
$$ E\left[\left\langle \mathbf{v}, \mathbf{P}_{kn} \right\rangle \right] \;=\; \sum\nolimits_{\sigma \in S_k} v_{\sigma} E\left[ P_{\sigma n} \right] \;=\; \sum\nolimits_{\sigma \in S_k} v_{\sigma} \tfrac{1}{k!} \;=\;  \left\langle \mathbf{v}, \mathbf{1} \right\rangle \tfrac{1}{k!} \;=\; 0
$$
Therefore $E[\left\langle \mathbf{v}, \mathbf{P}_{kn} \right\rangle^2] = V\left[\left\langle \mathbf{v}, \mathbf{P}_{kn} \right\rangle \right]$ as noted in the introduction.

Theorem~\ref{proj1} thus claims that the variance of the U-statistic $\left\langle \mathbf{v}, \mathbf{P}_{kn} \right\rangle$ has order of magnitude $n^{-r}$ for $\mathbf{v} \in V_r \setminus \{\mathbf{0}\}$. By Corollary~\ref{rank} it is sufficient to prove that the rank of this U-statistic is $r$.

\medskip
In order to to compute the rank, we investigate the functions $\{f_i\}_{0<i<k}$, defined in Section~\ref{back}. Substituting our kernel $f$ in that definition,
$$ f_i(x_1,\dots,x_i) \;=\; E\left[v_{\mathrm{perm}\left(x_1, \dots, x_i, X_{i+1}, \dots, X_k \right)}\right]$$
where $X_{i+1}, \dots, X_k$ are iid in $[0,1]^2$.

We will rewrite this function in a symmetrized form. 
Recall that we identify $S_{k-i}$ with the following subgroup of $S_k$: $$ S_{k-i} \;=\; \left\{\tau \in S_k \;\left|\; \tau(1)=1,\; \tau(2)=2,\;\dots,\; \tau(i)=i \right.\right\} $$
For such $\tau \in S_{k-i}$, we define $X_j^{\tau} := (Y_j,Z_{\tau(j)})$. Note that $(X_{i+1}, \dots, X_k)$ and $(X_{i+1}^{\tau}, \dots, X_k^{\tau})$ have the same joint distribution, because the $Y$s and $Z$s in each of these two sequences are $2(k-i)$ independent uniform random variables in the interval $[0,1]$. Therefore, for every $\tau \in S_{k-i}$,
$$ f_i(x_1,\dots,x_i) \;=\; E\left[v_{\mathrm{perm}\left(x_1, \dots, x_i, X_{i+1}^{\tau}, \dots, X_k^{\tau} \right)}\right] $$
This expression can be averaged over all $\tau \in S_{k-i}$,
$$ f_i(x_1,\dots,x_i) \;=\; \tfrac{1}{(k-i)!} \sum\limits_{\tau \in S_{k-i}} E\left[v_{\mathrm{perm}\left(x_1, \dots, x_i, X_{i+1}^{\tau}, \dots, X_k^{\tau} \right)}\right] $$
or, by linearity of expectation,
$$ f_i(x_1,\dots,x_i) \;=\; E\left[g_i\left(x_1, \dots, x_i, X_{i+1}, \dots, X_k \right)\right] $$
where
$$ g_i\left(x_1, \dots, x_k \right) \;:=\; \tfrac{1}{(k-i)!} \sum\limits_{\tau \in S_{k-i}} v_{\mathrm{perm}\left(x_1, \dots, x_i, x_{i+1}^{\tau}, \dots, x_k^{\tau} \right)} $$

Here are some observations about the function $g_i$. First, we show that for any generic input $x_1,\dots, x_k$, this function averages $v_{\sigma}$ over some two-sided coset of $S_{k-i}$ in $S_k$. Let $\alpha,\beta \in S_k$ encode the relative ordering of the coordinates of the points $x_j = (y_j,z_j)$ as follows: $y_{\beta(1)} < y_{\beta(2)} < \dots < y_{\beta(k)}$ and $z_{\alpha^{-1}(1)} < z_{\alpha^{-1}(2)} < \dots < z_{\alpha^{-1}(k)}$. In this notation,
$$ \mathrm{perm}\left(x_1, \dots, x_i, x_{i+1}^{\tau}, \dots, x_k^{\tau} \right) \;=\; \alpha \tau \beta $$
As $\tau \in S_{k-i}$ varies, we obtain all the $(k-i)!$ permutations that map $\beta^{-1}(j)$ to $\alpha(j)$ for $j\leq i$. In other words, $g_i\left(x_1, \dots, x_k \right)$ is the average of $v_{\sigma}$ over all $\sigma$ in the two-sided coset $\alpha S_{k-i}\beta$, where $\alpha$ and $\beta$ are determined by the input $x_1,\dots,x_k$.

Moreover, it is readily verified that the average $g_i(x_1,\dots,x_k)$ is symmetric under permuting $y_{i+1},\dots,y_k$ or permuting $z_{i+1},\dots,z_k$, as well as under permuting $x_1,\dots,x_i$. 

Finally, $g_i$ is invariant under order-preserving maps of the real numbers $y_1,\dots,y_k$ and under order-preserving maps of $z_1,\dots,z_k$, as it is defined in terms of $\mathrm{perm}(\cdots)$.

\medskip
The proof is completed in two step. We first prove that the rank of $\left\langle \mathbf{v}, \mathbf{P}_{kn} \right\rangle$ is at least $r$ by showing that $V[f_i] = 0$ for $i < r$. Then we will show that $V[f_r] > 0$, so that the rank is exactly $r$. 

\medskip
\noindent
\textbf{(I)}
Consider $f_i$ for $0<i<r$. The vector $\mathbf{v} \in V_r$ is a linear combination of matrix elements $\mathbf{R}^{\lambda}_{ij}$ for various $\lambda \vdash k$ with $\lambda_1 = k-r$. Letting $l = k-i > k-r$ puts $\mathbf{v}$ in the left hand side of Lemma~\ref{cosets}. The lemma implies that the sum of $v_{\sigma}$ over any two-sided coset of $S_l = S_{k-i}$ vanishes. 

As observed above, the function $g_i$ averages $v_{\sigma}$ over various two-sided cosets of $S_{k-i}$. We deduce that if $\mathbf{v} \in V_r$ for $r>i$ then $g_i(x_1,\dots, x_k) = 0$ on any generic input. 

In particular, $g_i(x_1,\dots,x_i,X_{i+1}, \dots, X_k)=0$ almost surely for generic $x_1,\dots,x_i$. Since $f_i$ is the expected value of $g_i$, also $f_i(x_1,\dots,x_i)=0$ on generic input. It follows that $f_i(X_1,\dots,X_i)=0$ almost surely, and $V[f_i]$ vanishes as desired.

\medskip
\noindent
\textbf{(II)}
We now consider the function $f_r$. The vector $\mathbf{v} \in V_r \setminus \{\mathbf{0}\}$ is \emph{not} contained in $V_{r+1} + V_{r+2} + \dots + V_{k-1}$, which is the left hand side of Lemma~\ref{cosets} with $l=k-r$. The lemma guarantees that the sum of $v_{\sigma}$ does \emph{not} always vanish simultaneously on all two-sided cosets of $S_{k-r}$.

Fix one such coset $C = \alpha S_{k-r} \beta$, where $\alpha, \beta \in S_k$, for which~$\sum_{\sigma \in C} v_{\sigma} \neq 0 $. Recall that the $r$ entries $\sigma(\beta^{-1}(i)) = \alpha(i)$ for $i \in \{1,\dots,r\}$ are common to all $\sigma \in C$. Denote the common $r$-pattern occurring in such $\sigma$ at these $r$ positions  by~$\rho \in S_r$. 

We will focus on the restriction of $f_r(x_1,\dots,x_r)$ to generic inputs that satisfy $\mathrm{perm}\left(x_1, \dots, x_r\right) = \rho$, where $\rho$ is determined by the choice of $C$ which depends on $\mathbf{v}$ as above. As usual, the order statistics of the input coordinates $x_j = (y_j,z_j) \in [0,1]^2$ are 
\begin{alignat*}{5}
0 \;&<\; y_{(1)} \;&<&\; y_{(2)} \;&<&\; \;\cdots\; \;&<&\; y_{(r)} \;&<&\; 1 \\
0 \;&<\; z_{(1)} \;&<&\; z_{(2)} \;&<&\; \;\cdots\; \;&<&\; z_{(r)} \;&<&\; 1 
\end{alignat*}
In these terms
\begin{align*}
& f_r\left( x_1, \dots, x_r \right) \;=\; E\left[g_r \left(x_1, \dots, x_r, X_{r+1}, \dots, X_k \right) \right] \\ 
& \;\;\;\;\;\;\;\; \;=\; E\left[ g_r\left((y_{(1)}, z_{(\rho(1))}),\dots,(y_{(r)}, z_{(\rho(r))}), (Y_{r+1},Z_{r+1}),\dots,(Y_k,Z_k)\right) \right]
\end{align*}
where $X_j = (Y_j,Z_j)$ are independent and uniform in $[0,1]^2$, and we permuted $x_1,\dots,x_r$ as observed above.

It also follows from the above observed symmetries of $g_r$ that the random variable $g_r (x_1, \dots, x_r, X_{r+1}, \dots, X_k )$ only depends on the two sets $\{Y_{r+1},\dots,Y_k\}$ and $\{Z_{r+1},\dots,Z_k\}$, regardless of how they are matched. Moreover, by the invariance to order-preserving maps, it only depends on how many of $\{Y_{r+1},\dots,Y_k\}$ lie in each interval between the fixed $y_{(1)},\dots,y_{(r)}$ and the same for the $z$-coordinate.

This means that $f_r$, the expected value of $g_r \left(x_1, \dots, x_r, X_{r+1}, \dots, X_k \right)$, can be computed by considering a finite number of options for how the $k-r$ random points are distributed between the $r+1$ possible $y$-intervals, and how between the $z$-intervals. Suppressing the variables $x_1, \dots, x_r$ we write
$$
f_r \;=\; \sum\limits_{\substack{i_0,i_1,\dots,i_r\geq 0 \\ \sum_t{i_t} = k-r}} \; \sum\limits_{\substack{j_0,j_1,\dots,j_r\geq 0 \\ \sum_t{j_t} = k-r}} \mathrm{Pr}\left[A_{i_0,\dots,i_r,j_0,\dots,j_r}\right] \, G_{i_0,\dots,i_r,j_0,\dots,j_r}
$$
where $A_{i_0,\dots,i_r,j_0,\dots,j_r}\left( x_1, \dots, x_r \right)$ is the event,
$$ \forall \; t \in \{0,\dots,r\} \;\;\; \left\{ \;\; \begin{matrix}
\left| \{Y_{r+1},\dots,Y_k\} \cap \left(y_{(t)},y_{(t+1)}\right) \right| \;=\; i_t \vspace{0.5em} \\
\left| \{Z_{r+1},\dots,Z_k\} \cap \left(z_{(t)},z_{(t+1)}\right) \right| \;=\; j_t
\end{matrix}\right.$$
and $G_{i_0,\dots,i_r,j_0,\dots,j_r}\left( x_1, \dots, x_r \right)$ is the common value of $g_r \left(x_1, \dots, x_r, X_{r+1}, \dots, X_k \right)$ given this event. Here we set $y_{(0)} = z_{(0)} = 0$ and $y_{(r+1)} = z_{(r+1)} = 1$. 

We go further and compute these probabilities. Instead of the given $r$ points, $x_1 = (y_1,z_1),\dots,x_r = (y_r,z_r)$, we define new variables corresponding to the differences between their reordered coordinates,
\begin{align*}
\Delta y_t \;&:=\; y_{(t+1)} - y_{(t)} \;=\; \mathrm{Pr}\left[y_{(t)} < Y < y_{(t+1)}\right] \\
\Delta z_t \;&:=\; z_{(t+1)} - z_{(t)} \;=\; \mathrm{Pr}\left[z_{(t)} < Z < z_{(t+1)}\right]
\end{align*}
where $(Y,Z)$ is uniform in $[0,1]^2$. Now the probabilities of the above events are monomial in these differences,
$$ \mathrm{Pr}\left[A_{i_0,\dots,i_r,j_0,\dots,j_r}\right] \;=\; \tbinom{k-r}{i_0 \, \dots \, i_r} (\Delta y_0)^{i_0} \cdots (\Delta y_r)^{i_r} \, \tbinom{k-r}{j_0 \, \dots \, j_r} (\Delta z_0)^{j_0} \cdots (\Delta z_r)^{j_r} $$
Here the multinomial coefficient $\tbinom{k-r}{i_0 \, \dots \, i_r}$ is the number of ways to assign $i_0,\dots,i_r$ of the $k-r$ random points to the $r+1$ respective $y$-intervals, and same for $j_0,\dots,j_r$ and $z$-intervals. Note that by definition no two of the events $A_{i_0,\dots,i_r,j_0,\dots,j_r}$ yield the same monomial in the variables $\Delta y_t$ and~$\Delta z_t$.

As for the other factor $G_{i_0,\dots,i_r,j_0,\dots,j_r}$, which is the corresponding value of~$g_r$, we have seen in the first part of the proof that this is the average of $v_{\sigma}$ over a two-sided coset of $S_{k-r}$. Namely,
$$ G_{i_0,\dots,i_r,j_0,\dots,j_r}(x_1,\dots,x_r) \;=\; \tfrac{1}{(k-r)!} \sum\limits_{\tau \in S_{k-r}} v_{\mathrm{perm}\left(x_1, \dots, x_r, X_{r+1}^{\tau}, \dots, X_k^{\tau} \right)} $$
for any $X_{r+1}, \dots, X_k$ that satisfy the event $A_{i_0,\dots,i_r,j_0,\dots,j_r}$. We next show that for some $i_0,\dots,i_r$ and $j_0,\dots,j_r$ it doesn't vanish.

Recall that we are interested in the domain where $\mathrm{perm}\left(x_1, \dots, x_r\right) = \rho$. This pattern $\rho \in S_r$ was determined by the two-sided coset $C = \alpha S_{k-r} \beta$, as the one that occurs in the positions $\beta^{-1}(1),\dots,\beta^{-1}(r)$ which are mapped by all members of~$C$ to $\alpha(1),\dots,\alpha(r)$ correspondingly. Let $i_t$ be the number of elements in $\{1,\dots,k\}$ that are larger than $t$ elements in $\{\beta^{-1}(1),\dots,\beta^{-1}(r)\}$ and smaller than the other $r-t$. Similarly, let $j_t$ be the number of elements in $\{1,\dots,k\}$ larger than $t$ elements in $\{\alpha(1),\dots,\alpha(r)\}$ and smaller than the other $r-t$. With this choice, any $x_{r+1}, \dots, x_k$ that satisfy the event $A_{i_0,\dots,i_r,j_0,\dots,j_r}$ fix $\beta^{-1}(1),\dots,\beta^{-1}(r)$ as the respective positions of $y_1,\dots,y_r$ in the sorted $\left\{y_1,\dots,y_k\right\}$, and similarly send $z_1,\dots,z_r$ to positions $\alpha(1),\dots,\alpha(r)$ respectively in the sorted $\left\{z_1,\dots,z_k\right\}$. Therefore, conditioning on this $A_{i_0,\dots,i_r,j_0,\dots,j_r}$, 
$$\left\{\mathrm{perm}\left(x_1, \dots, x_r, X_{r+1}^{\tau}, \dots, X_k^{\tau} \right)\right\}_{\tau \in S_{k-r}} \;=\; C $$ 
and hence
$$ G_{i_0,\dots,i_r,j_0,\dots,j_r}\left(x_1, \dots, x_r\right) \;=\; \tfrac{1}{(k-r)!} \sum\limits_{\sigma \in C} v_{\sigma} \;\neq\; 0 \;.$$

In conclusion, the above expression for $f_r$ is a polynomial of degree $2(k-r)$ in the $2(r+1)$ variables $\Delta y_0, \dots, \Delta y_r, \Delta z_0, \dots, \Delta z_r$. Moreover, at least one coefficient in this polynomial is nonzero. 

We claim that $f_r$ is also nonzero as a polynomial in the $2r$ original variables $y_1,\dots,y_r,z_1\dots,z_r$. Indeed this linear change of variables only restricts the domain to $\Delta y_0 + \dots + \Delta y_r = \Delta z_0 + \dots + \Delta z_r = 1$. But since our polynomial is both homogeneous in $\Delta y_0, \dots, \Delta y_r$ and homogeneous in $\Delta z_0, \dots, \Delta z_r$, this kind of restriction cannot make it vanish.

This nonzero polynomial form of $f_r\left((y_1, z_1), \dots, (y_r,z_r) \right)$ is only valid where $\mathrm{perm}\left(x_1, \dots, x_r\right) = \rho$ as determined by~$C$. Nevertheless, for uniform $X_1,\dots,X_r$ in $[0,1]^{2r}$ the event $\mathrm{perm}\left(X_1, \dots, X_r\right) = \rho$ is an open set that occurs with positive probability~$1/(r!)$. The polynomial $f_r$ is hence nonzero almost everywhere, conditioning on this event. 

We remark that, depending on~$\mathbf{v}$, the function $f_r$ might indeed vanish where its input induces certain $r$-patterns other than~$\rho$.
 
To sum up, $f_r\left(X_1, \dots, X_r\right)$ is nonzero with positive probability, while $E[f_r] = E[f] = 0$ since $\sum_{\sigma} v_\sigma = 0$ for $r>0$. This means that $V[f_r] > 0$, and the rank of $\left\langle \mathbf{v}, \mathbf{P}_{kn} \right\rangle$ is precisely~$r$.
\qed

\bigskip \noindent
\textbf{Proof of Theorem~\ref{cross}}
\nopagebreak \medskip \\ \noindent
Theorem~\ref{proj1} shows that the projection of the profile $\mathbf{P}_{kn}$ to $V_r$ has order of magnitude of~$n^{-r/2}$ along every direction in this subspace. Here we consider two normalized components $n^{r/2} \Pi_r \mathbf{P}_{kn}$ and $n^{s/2} \Pi_s \mathbf{P}_{kn}$ where $r<s$. A priori, these two random vectors in $V_r$ and $V_s$ are dependent, and in general correlated. The content of Theorem~\ref{cross} is that any such correlation tends to zero as~$n \to \infty$. 

It is sufficient to consider separately each entry of the cross-covariance matrix of $n^{r/2}\Pi_r \mathbf{P}_{kn}$ and $n^{s/2}\Pi_s \mathbf{P}_{kn}$. Since this matrix is the expected outer product of the two orthogonal projections, all we need to show is
$$ \mathrm{E} \left[\; n^{r/2} \left\langle \mathbf{u}, \mathbf{P}_{kn} \right\rangle \, n^{s/2} \left\langle \mathbf{v}, \mathbf{P}_{kn} \right\rangle \;\right] \; \xrightarrow{\: n \to \infty \: } \; 0 $$
for every $\mathbf{u} \in V_r$ and $\mathbf{v} \in V_s$. 
Recall that $\mathrm{E}\left[ \left\langle \mathbf{v}, \mathbf{P}_{kn} \right\rangle \right] = 0$ if $\mathbf{v} \in V_s$ for $s>0$, so that the above expectation of the product is indeed the covariance of the two projections.  

As before, we consider the profile of $\mathrm{perm}\left(X_1,\dots,X_n\right)$ where $X_1,\dots,X_n$ are independent and follow the continuous uniform distribution in $[0,1]^2$. Also, we similarly expand each scalar projection according to the $\tbinom{n}{k}$ occurrences of patterns,
\begin{align*}
\left\langle \mathbf{u}, \mathbf{P}_{kn} \right\rangle \left\langle \mathbf{v}, \mathbf{P}_{kn} \right\rangle
\;&=\;
\sum\limits_{\tau \in S_k} u_{\tau} P_{\tau} \sum\limits_{\sigma \in S_k} v_{\sigma} P_{\sigma} 
\\
\;&=\; \frac{1}{\tbinom{n}{k}^2} \sum\limits_{i_1 < \dots < i_k} \sum\limits_{j_1 < \dots < j_k} u_{\mathrm{perm}\left(X_{i_1},\dots,X_{i_k}\right)} \, v_{\mathrm{perm}\left(X_{j_1},\dots,X_{j_k}\right)} \;.
\end{align*}

The following argument applies to any product of two U-statistics of different ranks on the same sequence of iid random variables. We divide the terms of the double sum into two cases. Suppose that less than $s$ positions from $i_1,\dots,i_k$ and $j_1,\dots,j_k$ coincide. Then the expected value of this term can be written, up to ordering, as
\begin{align*}
&\mathrm{E} \left[ u_{\mathrm{perm}\left(X_{i_1},\dots,X_{i_k}\right)} \, v_{\mathrm{perm}\left(X_{i_1},\dots,X_{i_t},X_{j_{t+1}},\dots,X_{j_k}\right)}\right]
\\
&\;=\; \mathrm{E} \left[ u_{\mathrm{perm}\left(X_{i_1},\dots,X_{i_k}\right)} \, \mathrm{E} \left[ \left. v_{\mathrm{perm}\left(X_{i_1},\dots,X_{i_t},X_{j_{t+1}},\dots,X_{j_k}\right)}\;\right|\; X_{i_1},\dots,X_{i_k} \right]\right] \;.
\end{align*}
where $i_1,\dots,i_k$ and $j_{t+1},\dots,j_k$ are disjoint. By the proof of Theorem~\ref{proj1}, the rank of the kernel $v_{\mathrm{perm}(\cdot)}$ is $s$. Since $t<s$, for every generic $x_1,\dots,x_t$ and random variables $X_{t+1},\dots,X_k$,
$$ f_t(x_1,\dots,x_t) \;=\; \mathrm{E} \left[ v_{\mathrm{perm}\left(x_1,\dots,x_t,X_{t+1},\dots,X_k\right)}\right] \;=\; 0 \;.$$
because $V[f_t]=0$. Therefore the above conditional expectation is $0$ almost surely, and so is the expectation of the product. Thus terms of this kind vanish and do not contribute to the covariance.

In the remaining terms, at least $s$ positions from $i_1,\dots,i_k$ and $j_1,\dots,j_k$ coincide. We estimate the contribution of all these terms as, at most
$$ \binom{n}{s} \binom{n-s}{k-s} \binom{n-s}{k-s} \max\limits_{\tau} |u_{\tau}| \, \max\limits_{\sigma} |v_{\sigma}| \;=\; O\left(n^{2k-s}\right) \;.$$
Therefore, for $\mathbf{u} \in V_r$ and $\mathbf{v} \in V_s$ where $r<s$,
$$ \mathrm{E} \left[\; \left\langle \mathbf{u}, \mathbf{P}_{kn} \right\rangle \cdot \left\langle \mathbf{v}, \mathbf{P}_{kn} \right\rangle \;\right] \;=\; O\left(n^{-s}\right) \;=\; o\left(n^{-(r+s)/2}\right) $$  
and the theorem follows.
\qed

\bigskip \noindent
\textbf{Proof of Lemma~\ref{cosets}}
\nopagebreak \medskip \\ \noindent
We very briefly recall a construction of the simple representations of~$S_k$, which appears in Fulton and Harris~\cite[Lecture 4]{fulton1991representation}.

A \emph{Young tableau} $T$ of shape $\lambda \vdash k$ is a sequence of ordered \emph{rows} $T_1,\dots,T_\ell$ where $T_i = T_{i1}T_{i2} \dots T_{i\lambda_i}$ and $\bigcup_i T_i = \{1,\dots,k\}$. One \emph{Young subgroup} permutes the rows, $P_{T} = \{\sigma \in S_k\;|\; \sigma(T_{ij}) \in T_i\}$, and another one permutes the \emph{columns}, $Q_{T} = \{\sigma \in S_k\;|\; \sigma(T_{ij}) \in T_j'\}$ where $T'_j = T_{1j}T_{2j} \dots T_{\lambda_j'j}$ with $\lambda_j' = |\{i\,|\,j\leq\lambda_i\}|$. 

\begin{example*}
Let $T = (21,3)$. Then $P_T = \{123,213\}$, $T' = (23,1)$  and $Q_T = \{123,132\}$.
\end{example*}

Define a multiplication on the vector space $\mathbb{R}^{k!} = \{\sum_{\sigma \in S_k} v_{\sigma}\mathbf{\sigma}\;|\;v_{\sigma}\in\mathbb{R}\}$ extending $\mathbf{\sigma} \cdot \mathbf{\tau} = \mathbf{\sigma\tau}$. This turns it into an algebra $A = \mathbb{R}S_k$. 

Given a tableau $T$ of shape~$\lambda$, we consider two elements, $a_T = \sum_{\sigma \in P_T}\sigma$ and $b_T = \sum_{\sigma \in Q_T}\mathrm{sign}(\sigma)\sigma$, written also as $a_{\lambda}$ and $b_{\lambda}$ when the choice of $T$ is not important. The \emph{Young symmetrizer} is $c_{T} = c_{\lambda} = a_{\lambda} \cdot b_{\lambda}$. 

\begin{example*}
For $T$ as above, $c_T = (123+213) \cdot (123-132) = 123-132+213-231$.
\end{example*}

The simple representation $\rho^{\lambda}$ is now given by the left action of $S_k \subset A$ on the subspace spanned by the corresponding Young symmetrizer: $Ac_{\lambda} = \{x\cdot c_{\lambda}\,|\,x \in A\}$ \cite[Theorem 4.3]{fulton1991representation}. We use this construction to prove both directions of Lemma~\ref{cosets}.

\medskip
\noindent
\textbf{($\Rightarrow$)}
It is sufficient to show that every matrix element $\mathbf{R}_{ij}^{\lambda}$, for $\lambda \vdash k$ with $\lambda_1 < l$ and $1 \leq i,j \leq d_{\lambda}$, is contained in the right hand side of the lemma.

Let the partition $\mu \vdash k$ be $(l,1,1,\dots,1)$, and let the tableau $T$ have shape~$\mu$, with $T_1 = \{k,k-1,\dots,k-l+1\}$. Note that $P_{\mu} = S_l$ as defined before the lemma. Since $\lambda_1 < l = \mu_1$, it follows that $\lambda < \mu$ in the lexicographic order~\cite[4.22]{fulton1991representation}. For $\alpha,\beta \in S_k$ let
$$ f_{\alpha\beta} \;:=\; \sum_{\sigma \in \alpha S_l \beta} \mathbf{\sigma} \;=\; \mathbf{\alpha} \cdot \sum_{\tau \in P_{\mu}} \mathbf{\tau} \cdot \mathbf{\beta} \;=\; \mathbf{\alpha} \cdot a_{\mu} \cdot \mathbf{\beta} \;\in\; A $$
Then for every $x \in A$
$$ f_{\alpha\beta}\cdot x\cdot c_\lambda \;=\; \mathbf{\alpha} \cdot a_{\mu} \cdot (\mathbf{\beta} \cdot x \cdot a_\lambda) \cdot b_\lambda \;=\; \mathbf{\alpha} \cdot 0 \;=\; 0 $$ 
by~\cite[Lemma 4.23]{fulton1991representation} since $\lambda < \mu$. This means that $f_{\alpha\beta}$ annihilates all of~$Ac_{\lambda}$, the $A$-module that defines the representation $\rho^{\lambda}$ of $S_k$. For any choice of basis of~$Ac_{\lambda}$, the corresponding matrices $R^{\lambda} (\sigma) \in \mathbb{R}^{d_{\lambda}\times d_{\lambda}}$ must hence satisfy
$$ \forall \alpha,\beta \in S_k \;\;\;\;\;\;\;\;\;\; \sum_{\sigma \in \alpha S_l \beta} R^{\lambda} (\sigma) \;=\; \mathbf{0} $$
In other words, the average of every matrix entry $R_{ij}^{\lambda}(\sigma)$ with $\lambda_1 < l$, over~$\sigma$ in any two-sided coset of $S_l$, vanishes as required.

\medskip
\noindent
\textbf{($\Leftarrow$)}
Let $T$ be a tableau with shape $\lambda \vdash k$ with $\lambda_1 \geq l$. The subgroup~$P_{\lambda}$ includes a smaller subgroup, all permutations that only shuffle the first $l$ entries of the first row~$T_1$. This subgroup and each of its cosets are two-sided cosets of $S_l$. Therefore $P_{\lambda} = \mathaccent\cdot\cup_{(\alpha,\beta)\in\mathcal{I}}\, \alpha S_l \beta$ for a suitable $\mathcal{I} \subset S_k \times S_k$.

Let $v = \sum_{\sigma \in S_k}v_{\sigma}\mathbf{\sigma} \in A$ and assume that $\sum v_{\sigma}$ vanishes over every two-sided coset of $S_l$, as in the right hand side of the lemma. For every $\gamma \in S_k$,
$$ v \cdot \mathbf{\gamma} \cdot c_{\lambda} \;=\; \sum\limits_{\sigma \in S_k} v_{\sigma} \mathbf{\sigma \gamma} \cdot \sum\limits_{(\alpha,\beta)\in\mathcal{I}} \;\sum\limits_{\tau \in S_l} \mathbf{\alpha \tau \beta} \cdot b_{\lambda} \;=\; $$
letting $\rho = \sigma\gamma\alpha\tau$,
$$ \;=\; \sum\limits_{\alpha,\beta} \;\sum\limits_{\tau \in S_l} \;\sum\limits_{\rho \in S_k} v_{(\rho(\gamma\alpha\tau)^{-1})} \mathbf{\rho \beta} \cdot b_{\lambda} \;=\; \sum\limits_{\alpha,\beta} \;\sum\limits_{\rho \in S_k} \left[\sum\limits_{\omega \in \rho S_l (\gamma\alpha)^{-1}} v_{\omega} \right]\mathbf{\rho \beta} \cdot b_{\lambda} \;=\; 0 $$
using the assumption on~$v$ for these $S_l$-cosets. By linearity $v$ annihilates $Ac_{\lambda}$, so that $\sum_{\sigma \in S_k} v_{\sigma} R^{\lambda}(\sigma) = \mathbf{0}$ similarly to the above. In other words $\langle \mathbf{v}, \mathbf{R}_{ij}^{\lambda} \rangle = 0$ for every $\lambda \vdash k$ with $\lambda_1 \geq l$ and $1 \leq i,j \leq d_{\lambda}$. By the orthogonality of the matrix elements of different representations, $\mathbf{v}$ is spanned by the other $\mathbf{R}_{ij}^{\lambda}$, with $\lambda_1 < l$, as desired.
\qed

\section{Diagonalization}
\label{diag}

We have seen that the random profile $\mathbf{P}_{kn}$ decomposes according to the orthogonal subspaces $V_0,\dots,V_{k-1}$, and its projection on $V_r$ has order of magnitude $n^{-r/2}$. We therefore define the \emph{normalized $k$-profile}: 
$$ \mathbf{\tilde{P}}_{kn} \;:=\; \sum\limits_{r=0}^{k-1} n^{r/2} \,\Pi_r \mathbf{P}_{kn} $$

In order to study the asymptotic structure of the normalized profile, we ask how its one-dimensional projections are related to each other. This is reasonably captured by the $k!$-by-$k!$ limit second-moment matrix, 
$$ M_{k} \;:=\; \lim\limits_{n \to \infty} \mathrm{E} \left[ \mathbf{\tilde{P}}_{kn}^{} \mathbf{\tilde{P}}_{kn}^{\;T}  \right] $$

In terms of $M_k$, Theorems~\ref{proj1} and \ref{cross} claim that $\mathbf{v}^T M_k \mathbf{v} > 0$ for every $\mathbf{v} \in \bigcup_r V_r$, and $\mathbf{u}^T M_k \mathbf{v} = 0$ for every $\mathbf{u} \in V_r$ and $\mathbf{v} \in V_s$ such that $r \neq s$. It follows that this matrix is positive definite. Also, $M_k$ is equivalent as a quadratic form to a block matrix with $k$ blocks, by a change of basis to vectors from these subspaces. In this form, $M_k$ and the limit of $\text{cov}\,[ \mathbf{\tilde{P}}_{kn} ]$ are the same except at the 1-by-1 block that corresponds to $V_0$, as noted after Theorem~\ref{proj1}.

As usual in multivariate analysis, the distribution of $\mathbf{\tilde{P}}_{kn}$ is better understood by means of an orthonormal basis that fully diagonalizes $M_k$. This requires a basis of directions within each $V_r$, along which we project the profile. This would thus refine the profile decomposition from $k$ to $k!$ components that are asymptotically uncorrelated. The eigenvalues would tell how the variance varies between different projections, and the eigenvectors may point at particular combinations of pattern densities that are interesting to look~at. 

Since $V_r$ is spanned by matrix elements of $S_k$-representations, we ask if there exist specific representatives in their similarity classes that yield such bases. The following proposition indicates that this is indeed the case. 
\begin{prop}\label{diagonal}
Let $k \leq 6$. There exist an orthogonal $k!$-by-$k!$ matrix $U_k$ and a~strictly positive diagonal matrix $D_k$ such that
$$ U_k^T M_k U_k \;=\; D_k $$
and the columns of $U_k$ are the normalized matrix elements, 
$$\left\{\sqrt{\tfrac{d_{\lambda}}{k!}}\;\mathbf{R}^{\lambda}_{ij}\right\}_{\lambda \vdash k,\;1 \leq i,j \leq d_{\lambda}}$$
of orthogonal simple $S_k$-representations $\{R^\lambda\}_{\lambda \vdash k}$.  
\end{prop}
Here the matrices $R^{\lambda}(\sigma)$ for $\sigma \in S_k$ are orthogonal with respect to the standard inner products in $\mathbb{R}^{d_{\lambda}}$, and therefore the matrix $U_k$ is orthogonal with respect to the standard inner product in $\mathbb{R}^{k!}$~\cite[2.35]{fulton1991representation}. 

For finite $n$, Proposition~\ref{diagonal} suggests using these $k!$ scalar projections to express the $k$-profile:
$$ U_k^T\,\mathbf{\tilde{P}}_{kn} \;=\;\left(\;\left\langle n^{(k-\lambda_1)/2} \sqrt{\tfrac{d_{\lambda}}{k!}}\;\mathbf{R}^{\lambda}_{ij}\,,\,\mathbf{P}_{kn}\right\rangle\;\right)_{\lambda \vdash k,\; 1 \leq i,j \leq d_\lambda} $$

\begin{proof}
The proof of Proposition~\ref{diagonal} was found by computer-aided exploration in \texttt{Python}, and verified by the computer algebra system \texttt{Sage}~\cite{sagemath}. See \cite[\texttt{verify.sage}]{even2018gh} for the implementation. The file \texttt{main.sage} was run to perform the verification, and produced the output given in \texttt{output.txt}. The computation is straightforward, and described below.

We start with the second-moment matrix of the unnormalized profile. Consider the matrix $\tbinom{n}{k} \,\mathrm{E} [ \mathbf{P}_{kn}^{} \mathbf{P}_{kn}^{\;T}]$. Each entry of this matrix is a polynomial in~$\mathbb{Q}[n]$ of degree at most~$k$, as one can see in the proof of Theorem~\ref{cross}. We extrapolate these polynomials from the points $n = k,\dots,2k$, for which the second-moment matrices are computed directly by averaging the outer product over~$S_n$. 

Then we construct $U_k$. For every partition $\lambda \vdash k$, two matrices are given: $R^{\lambda}(\tau_k), R^{\lambda}(\rho_k) \in \overline{\mathbb{Q}}\,^{d_\lambda \times d_\lambda}$, where $ \tau_k = 2134...k$ and $\rho_k = 234...k1$. Since these transposition and rotation generate $S_k$, we can find $R^{\lambda}(\sigma)$ for all \mbox{$\sigma \in S_k$}. Our choice of generators of every simple representation $R^{\lambda}$, for $\lambda \vdash k \leq 6$, appears in the file \texttt{generators.sage}. All the entries in these matrices have finite descriptions as field elements in quadratic extensions of~$\mathbb{Q}$. The matrix elements of all  simple representations of~$S_k$ yield the orthogonal matrix~$U_k$ as described in the proposition. 

Then we compute the conjugate matrix $U^T_k \mathrm{E} [ \mathbf{P}_{kn}^{} \mathbf{P}_{kn}^{\;T}] U_k$. The leading order in~$n$ of every entry in this matrix is bounded according to the block structure, as follows from Theorems~\ref{proj1} and~\ref{cross}. Specifically, if $\mathbf{u}$ and $\mathbf{v}$ are matrix elements from $V_r$ and~$V_s$ respectively, then the entry $\mathbf{u}^T \mathrm{E} [ \mathbf{P}_{kn}^{} \mathbf{P}_{kn}^{\;T}] \mathbf{v}$ is $O(n^{-\max\left(r,s\right)})$. 

We thus obtain the normalized matrix $U_k^T \mathrm{E} [ \mathbf{\tilde{P}}_{kn}^{} \mathbf{\tilde{P}}_{kn}^{\;T}  ] U_k$ by multiplying the different blocks by appropriate powers of $\sqrt{n}$. The entries of this matrix are rational functions in~$\overline{\mathbb{Q}}(\sqrt{n})$. Taking $n \to \infty$ yields the leading, constant term, which is $U_k^TM_kU_k$. All that remains is to verify that this matrix is diagonal.
\end{proof}

Proposition~\ref{diagonal} and its proof raise new directions for further investigation, and some of them are discussed in the following examples and remarks.

\begin{example*}
We demonstrate the verification procedure for $k=2$. A direct computation yields
$$
\mathrm{E} [ \mathbf{P}_{22}^{} \mathbf{P}_{22}^{\;T}]=
\left(
\begin{smallmatrix}
\tfrac{1}{2} & 0 \\
0 & \tfrac{1}{2}
\end{smallmatrix}
\right)
,\;\;\;
\mathrm{E} [ \mathbf{P}_{23}^{} \mathbf{P}_{23}^{\;T}]=
\left(
\begin{smallmatrix}
\tfrac{19}{54} & \tfrac{4}{27} \\
\tfrac{4}{27} & \tfrac{19}{54}
\end{smallmatrix}
\right),\;\;\;
\mathrm{E} [ \mathbf{P}_{24}^{} \mathbf{P}_{24}^{\;T}]=
\left(
\begin{smallmatrix}
\tfrac{67}{216} & \tfrac{41}{216} \\
\tfrac{41}{216} & \tfrac{67}{216}
\end{smallmatrix}
\right)
$$
and by extrapolation,
$$
\binom{n}{2}
\mathrm{E} [ \mathbf{P}_{2n}^{} \mathbf{P}_{2n}^{\;T}]
\;=\;
\left(
\begin{smallmatrix}
\tfrac18 n^2 - \tfrac{5}{72} n + \tfrac{5}{36} \;&\;
\tfrac18 n^2 - \tfrac{13}{72} n - \tfrac{5}{36} \\[0.3em]  
\tfrac18 n^2 - \tfrac{13}{72} n - \tfrac{5}{36} \;&\; 
\tfrac18 n^2 - \tfrac{5}{72} n + \tfrac{5}{36} 
\end{smallmatrix}
\right)\;.
$$
The matrix elements of the trivial and the alternating representations of $S_2$ give $U_2 = \tfrac{1}{\sqrt{2}} \left( \begin{smallmatrix}
1 & 1
\\[0.3em]
1 & -1
\end{smallmatrix} \right)$, and the normalization is $U_2^T \mathbf{\tilde{P}}_{2n} =
\left(
\begin{smallmatrix}
1 & 0
\\[0.3em]
0 & \sqrt{n}
\end{smallmatrix}
\right) U_2^T \mathbf{P}_{2n}$. 
It follows that 
$$ U^T_2 \mathrm{E} [ \mathbf{\tilde{P}}_{2n}^{} \mathbf{\tilde{P}}_{2n}^{\;T}] U_2 = \frac{1}{\tbinom{n}{2}}\left(
\begin{array}{cc}
\tfrac{1}{4} n^2 - \tfrac{1}{4} n &  
0
\\[0.3em]
0 &  
\tfrac{1}{9} n^2 + \tfrac{5}{18} n
\end{array}
\right) \xrightarrow{\: n \to \infty \: } \; \left(\begin{array}{cc}
\tfrac{1}{2} & 0
\\[0.3em]
0 & \tfrac{2}{9}
\end{array}\right) $$
which is diagonal of full rank as required. This simple case is not new because all blocks in $U^T_2 M_2 U_2$ are one-by-one, but in general the diagonalization is verified also within blocks.
\end{example*}

\begin{example*}
In the case $k=3$ the choice of representation matrices becomes important. For the standard representation corresponding to $\lambda = (2,1)$ we use:
\def\arraystretch{1.5}
\begin{align*}
&
R^{21}(\tau_{3}) = 
\left(\begin{array}{cc}
-\tfrac{1}{2} & \tfrac{\sqrt{3}}{2} \\
\tfrac{\sqrt{3}}{2} & \tfrac{1}{2}
\end{array}\right)
&&
R^{21}(\rho_{3}) = 
\left(\begin{array}{cc}
-\tfrac{1}{2} & -\tfrac{\sqrt{3}}{2} \\
\tfrac{\sqrt{3}}{2} & -\tfrac{1}{2}
\end{array}\right)
\end{align*} 
As always, the trivial representation is given by $R^{k}(\tau_k) = R^{k}(\rho_k) = [1]$, and the alternating by $R^{1\cdots 1}(\tau_k) = [-1]$  and $R^{1 \cdots 1}(\rho_k) = [(-1)^{k-1}]$. Repeating the above procedure, the diagonalized normalized 3-profile is
\def\arraystretch{1.35}
$$ U_3^T \mathbf{\tilde{P}}_{3n} \;=\; \left( 
\begin{smallmatrix}
\tfrac{1}{\sqrt{6}} & \tfrac{\sqrt{n}}{\sqrt{3}} & 0 & 0 & \tfrac{\sqrt{n}}{\sqrt{3}} & \tfrac{n}{\sqrt{6}} \\
\tfrac{1}{\sqrt{6}} & -\tfrac{\sqrt{n}}{2  \sqrt{3}} & \,-\tfrac{\sqrt{n}}{2}\, & \,-\tfrac{\sqrt{n}}{2}\, & \tfrac{\sqrt{n}}{2  \sqrt{3}} & -\tfrac{n}{\sqrt{6}} \\
\tfrac{1}{\sqrt{6}} & -\tfrac{\sqrt{n}}{2  \sqrt{3}} & \tfrac{\sqrt{n}}{2} & \tfrac{\sqrt{n}}{2} & \tfrac{\sqrt{n}}{2  \sqrt{3}} & -\tfrac{n}{\sqrt{6}} \\
\tfrac{1}{\sqrt{6}} & -\tfrac{\sqrt{n}}{2  \sqrt{3}} & -\tfrac{\sqrt{n}}{2} & \tfrac{\sqrt{n}}{2} & -\tfrac{\sqrt{n}}{2  \sqrt{3}} & \tfrac{n}{\sqrt{6}} \\
\tfrac{1}{\sqrt{6}} & -\tfrac{\sqrt{n}}{2  \sqrt{3}} & \tfrac{\sqrt{n}}{2} & -\tfrac{\sqrt{n}}{2} & -\tfrac{\sqrt{n}}{2  \sqrt{3}} & \tfrac{n}{\sqrt{6}} \\
\tfrac{1}{\sqrt{6}} & \tfrac{\sqrt{n}}{\sqrt{3}} & 0 & 0 & -\tfrac{\sqrt{n}}{\sqrt{3}} & -\tfrac{n}{\sqrt{6}}
\end{smallmatrix}
\right)^T \left( 
\begin{array}{c}
P_{123,n} \\
P_{132,n} \\
P_{213,n} \\
P_{231,n} \\
P_{312,n} \\
P_{321,n}
\end{array} \right) $$
The verification proceeds by computing its second moment matrix $U^T_3 \mathrm{E} [ \mathbf{\tilde{P}}_{3n}^{} \mathbf{\tilde{P}}_{3n}^{\;T}] U_3$ as before, and confirming that as $n \to \infty$ its limit $U^T_3 M_3 U_3$ is diagonal of full rank. 
\end{example*}

\begin{example*}
Our choice of generators for $\lambda \vdash 4$ is given below. For $\lambda \vdash 5$ see Table~\ref{gen5}. For all $\lambda \vdash k \leq 6$ see~\cite[\texttt{generators.sage}]{even2018gh}.

\def\arraystretch{1.25}
\begin{align*}
&
R^{31}(\tau_{4}) = 
\left(\begin{array}{ccc}
\tfrac{1}{5} & \tfrac{2}{\sqrt{5}} & -\tfrac{2}{5} \\
\tfrac{2}{\sqrt{5}} & 0 & \tfrac{1}{\sqrt{5}} \\
-\tfrac{2}{5} & \tfrac{1}{\sqrt{5}} & \tfrac{4}{5}
\end{array}\right)
&&
R^{31}(\rho_{4}) = 
\left(\begin{array}{ccc}
-\tfrac{4}{5} & -\tfrac{1}{\sqrt{5}} & -\tfrac{2}{5} \\
\tfrac{1}{\sqrt{5}} & 0 & -\tfrac{2}{\sqrt{5}} \\
-\tfrac{2}{5} & \tfrac{2}{\sqrt{5}} & -\tfrac{1}{5}
\end{array}\right)
\\[0.5em] &
R^{22}(\tau_{4}) = 
\left(\begin{array}{cc}
- 1 & 0 \\
0 & 1
\end{array}\right)
&&
R^{22}(\rho_{4}) = 
\left(\begin{array}{cc}
\tfrac{1}{2} & \tfrac{\sqrt{3}}{2} \\
\tfrac{\sqrt{3}}{2} & -\tfrac{1}{2}
\end{array}\right)
\\[0.5em] &
R^{211}(\tau_{4}) = 
\left(\begin{array}{ccc}
0 & 0 & 1 \\
0 & - 1 & 0 \\
1 & 0 & 0
\end{array}\right)
&&
R^{211}(\rho_{4}) = 
\left(\begin{array}{ccc}
0 & 1 & 0 \\
- 1 & 0 & 0 \\
0 & 0 & 1
\end{array}\right)
\end{align*}
\end{example*}

\begin{table}[p]
\def\arraystretch{1.5625}
\setlength\arraycolsep{3pt}
\noindent
\begin{tabular}{|c||c|c|}
\hline
$\lambda$ & $R^{\lambda}(\tau_5)$ & $R^{\lambda}(\rho_5)$
\\
\hline
\hline
5 & 1 & 1 \\
\hline
$\begin{array}{c} 4 \\ 1 \end{array}$
&
$\begin{array}{cccc}
\tfrac{1}{10} & -\tfrac{3}{10} & \tfrac{3}{2  \sqrt{7}} & \tfrac{9}{2  \sqrt{35}} \\
-\tfrac{3}{10} & \tfrac{9}{10} & \tfrac{1}{2  \sqrt{7}} & \tfrac{3}{2  \sqrt{35}} \\
\tfrac{3}{2  \sqrt{7}} & \tfrac{1}{2  \sqrt{7}} & \tfrac{9}{14} & -\tfrac{3  \sqrt{5}}{14} \\
\tfrac{9}{2  \sqrt{35}} & \tfrac{3}{2  \sqrt{35}} & -\tfrac{3  \sqrt{5}}{14} & \tfrac{5}{14}
\end{array}$
&
$\begin{array}{cccc}
-\tfrac{1}{2} & -\tfrac{1}{2} & \tfrac{3}{2  \sqrt{7}} & -\tfrac{\sqrt{5}}{2  \sqrt{7}} \\
-\tfrac{1}{2} & 0 & \tfrac{1}{2  \sqrt{7}} & \tfrac{\sqrt{5}}{\sqrt{7}} \\
-\tfrac{3}{2  \sqrt{7}} & -\tfrac{1}{2  \sqrt{7}} & -\tfrac{11}{14} & -\tfrac{\sqrt{5}}{14} \\
\tfrac{\sqrt{5}}{2  \sqrt{7}} & -\tfrac{\sqrt{5}}{\sqrt{7}} & -\tfrac{\sqrt{5}}{14} & \tfrac{2}{7}
\end{array}$
\\  & & \\[-1.5em]\hline
$\begin{array}{c} 3 \\ 2 \end{array}$
&
$\begin{array}{ccccc}
-\tfrac{6}{7} & 0 & -\tfrac{1}{\sqrt{7}} & 0 & -\tfrac{\sqrt{6}}{7} \\
0 & 1 & 0 & 0 & 0 \\
-\tfrac{1}{\sqrt{7}} & 0 & 0 & 0 & \tfrac{\sqrt{6}}{\sqrt{7}} \\
0 & 0 & 0 & 1 & 0 \\
-\tfrac{\sqrt{6}}{7} & 0 & \tfrac{\sqrt{6}}{\sqrt{7}} & 0 & -\tfrac{1}{7}
\end{array}$
&
$\begin{array}{ccccc}
\tfrac{2}{7} & \tfrac{3}{2  \sqrt{7}} & -\tfrac{2}{\sqrt{7}} & -\tfrac{1}{2  \sqrt{14}} & -\tfrac{\sqrt{3}}{14  \sqrt{2}} \\
-\tfrac{3}{2  \sqrt{7}} & -\tfrac{1}{2} & -\tfrac{1}{2} & -\tfrac{1}{2  \sqrt{2}} & -\tfrac{\sqrt{3}}{2  \sqrt{14}} \\
\tfrac{2}{\sqrt{7}} & -\tfrac{1}{2} & 0 & -\tfrac{1}{2  \sqrt{2}} & -\tfrac{\sqrt{3}}{2  \sqrt{14}} \\
-\tfrac{1}{2  \sqrt{14}} & \tfrac{1}{2  \sqrt{2}} & \tfrac{1}{2  \sqrt{2}} & -\tfrac{1}{4} & -\tfrac{5  \sqrt{3}}{4  \sqrt{7}} \\
-\tfrac{\sqrt{3}}{14  \sqrt{2}} & \tfrac{\sqrt{3}}{2  \sqrt{14}} & \tfrac{\sqrt{3}}{2  \sqrt{14}} & -\tfrac{5  \sqrt{3}}{4  \sqrt{7}} & \tfrac{13}{28}
\end{array}$
\\ & & \\[-1.5em] \hline
$\begin{array}{c} 3 \\ 1 \\ 1 \end{array}$
&
\scalebox{0.875}{$\begin{array}{cccccc}
\tfrac{2}{5} & -\tfrac{\sqrt{6}}{5} & \tfrac{\sqrt{3}}{\sqrt{5}} & 0 & 0 & 0 \\
-\tfrac{\sqrt{6}}{5} & \tfrac{3}{5} & \tfrac{\sqrt{2}}{\sqrt{5}} & 0 & 0 & 0 \\
\tfrac{\sqrt{3}}{\sqrt{5}} & \tfrac{\sqrt{2}}{\sqrt{5}} & 0 & 0 & 0 & 0 \\
0 & 0 & 0 & 0 & \tfrac{1}{\sqrt{7}} & -\tfrac{\sqrt{6}}{\sqrt{7}} \\
0 & 0 & 0 & \tfrac{1}{\sqrt{7}} & -\tfrac{6}{7} & -\tfrac{\sqrt{6}}{7} \\
0 & 0 & 0 & -\tfrac{\sqrt{6}}{\sqrt{7}} & -\tfrac{\sqrt{6}}{7} & -\tfrac{1}{7}
\end{array}$}
&
\scalebox{0.875}{$\begin{array}{cccccc}
1 & 0 & 0 & 0 & 0 & 0 \\
0 & -\tfrac{1}{4} & -\tfrac{\sqrt{5}}{6  \sqrt{2}} & -\tfrac{5}{6  \sqrt{2}} & \tfrac{5}{2  \sqrt{14}} & \tfrac{5}{4  \sqrt{21}} \\
0 & \tfrac{\sqrt{5}}{6  \sqrt{2}} & \tfrac{1}{6} & \tfrac{\sqrt{5}}{6} & \tfrac{\sqrt{5}}{6  \sqrt{7}} & \tfrac{5  \sqrt{5}}{2  \sqrt{42}} \\
0 & \tfrac{5}{6  \sqrt{2}} & \tfrac{\sqrt{5}}{6} & -\tfrac{2}{3} & -\tfrac{2}{3  \sqrt{7}} & \tfrac{1}{2  \sqrt{42}} \\
0 & \tfrac{5}{2  \sqrt{14}} & -\tfrac{\sqrt{5}}{6  \sqrt{7}} & \tfrac{2}{3  \sqrt{7}} & \tfrac{4}{7} & -\tfrac{13}{14  \sqrt{6}} \\
0 & \tfrac{5}{4  \sqrt{21}} & -\tfrac{5  \sqrt{5}}{2  \sqrt{42}} & -\tfrac{1}{2  \sqrt{42}} & -\tfrac{13}{14  \sqrt{6}} & \tfrac{5}{28}
\end{array}$}
\\  & & \\[-1.5em]\hline
$\begin{array}{c} 2 \\ 2 \\ 1 \end{array}$
&
$\begin{array}{ccccc}
-\tfrac{1}{5} & \tfrac{2}{5} & 0 & \tfrac{2}{\sqrt{15}} & \tfrac{2  \sqrt{2}}{\sqrt{15}} \\
\tfrac{2}{5} & -\tfrac{4}{5} & 0 & \tfrac{1}{\sqrt{15}} & \tfrac{\sqrt{2}}{\sqrt{15}} \\
0 & 0 & 1 & 0 & 0 \\
\tfrac{2}{\sqrt{15}} & \tfrac{1}{\sqrt{15}} & 0 & -\tfrac{2}{3} & \tfrac{\sqrt{2}}{3} \\
\tfrac{2  \sqrt{2}}{\sqrt{15}} & \tfrac{\sqrt{2}}{\sqrt{15}} & 0 & \tfrac{\sqrt{2}}{3} & -\tfrac{1}{3}
\end{array}$
&
$\begin{array}{ccccc}
\tfrac{3}{10} & -\tfrac{1}{10} & \tfrac{3  \sqrt{3}}{2  \sqrt{10}} & \tfrac{\sqrt{3}}{2  \sqrt{5}} & \tfrac{\sqrt{3}}{2  \sqrt{10}} \\
-\tfrac{1}{10} & -\tfrac{4}{5} & -\tfrac{\sqrt{3}}{2  \sqrt{10}} & \tfrac{2}{\sqrt{15}} & -\tfrac{1}{2  \sqrt{30}} \\
-\tfrac{3  \sqrt{3}}{2  \sqrt{10}} & \tfrac{\sqrt{3}}{2  \sqrt{10}} & \tfrac{1}{4} & \tfrac{1}{2  \sqrt{2}} & -\tfrac{1}{4} \\
-\tfrac{\sqrt{3}}{2  \sqrt{5}} & -\tfrac{2}{\sqrt{15}} & \tfrac{1}{2  \sqrt{2}} & -\tfrac{2}{3} & \tfrac{1}{6  \sqrt{2}} \\
-\tfrac{\sqrt{3}}{2  \sqrt{10}} & \tfrac{1}{2  \sqrt{30}} & -\tfrac{1}{4} & \tfrac{1}{6  \sqrt{2}} & \tfrac{11}{12}
\end{array}$
\\ & & \\[-1.5em] \hline
$\begin{array}{c} 2 \\ 1 \\ 1 \\ 1 \end{array}$
&
$\begin{array}{cccc}
-\tfrac{1}{6} & \tfrac{\sqrt{5}}{2  \sqrt{3}} & \tfrac{\sqrt{5}}{6} & \tfrac{\sqrt{5}}{2  \sqrt{3}} \\
\tfrac{\sqrt{5}}{2  \sqrt{3}} & -\tfrac{1}{2} & \tfrac{1}{2  \sqrt{3}} & \tfrac{1}{2} \\
\tfrac{\sqrt{5}}{6} & \tfrac{1}{2  \sqrt{3}} & -\tfrac{5}{6} & \tfrac{1}{2  \sqrt{3}} \\
\tfrac{\sqrt{5}}{2  \sqrt{3}} & \tfrac{1}{2} & \tfrac{1}{2  \sqrt{3}} & -\tfrac{1}{2}
\end{array}$
&
$\begin{array}{cccc}
-\tfrac{2}{3} & 0 & \tfrac{\sqrt{5}}{6} & \tfrac{\sqrt{5}}{2  \sqrt{3}} \\
0 & 0 & \tfrac{\sqrt{3}}{2} & -\tfrac{1}{2} \\
\tfrac{\sqrt{5}}{6} & -\tfrac{\sqrt{3}}{2} & \tfrac{1}{6} & \tfrac{1}{2  \sqrt{3}} \\
-\tfrac{\sqrt{5}}{2  \sqrt{3}} & -\tfrac{1}{2} & -\tfrac{1}{2  \sqrt{3}} & -\tfrac{1}{2}
\end{array}$
\\ & & \\[-1.5em] \hline
$1^5$ & $-1$ & $1$
\\ \hline
\end{tabular}
~ \\ ~
\caption{Generators for representations of $S_5$ that verify Proposition~\ref{diagonal}.}
\label{gen5}
\end{table}

\begin{rmk*}
Janson, Nakamura and Zeilberger note that the eigenvalues of the first, $(k-1)^2$-dimensional component are rational numbers, for $k \leq 5$ \cite[Example 4.8]{janson2015asymptotic}. Here we report that this property persists across the whole $k!$-point spectrum of the normalized profile, and for $k \leq 6$~\cite[\texttt{verify-eigenvalues.sage}]{even2018gh}.

The eigenvectors, which are the matrix elements of the representations that we found, are only populated by signed square roots of rational numbers. In fact, for every $k \leq 6$ these are contained in the quadratic extension 
$$\mathbb{Q}_k \;:=\; \mathbb{Q}\left(\left\{\sqrt{p} \;|\; 2 \leq p < 2(k-1)\right\}\right)\;.$$

More structural properties of these representation matrices were found, and also used, in the process that led to their discovery. We do not discuss them further in this paper.
\end{rmk*}

\begin{rmk*}
The representation matrices that we use are not those of the known constructions: Young's semi-normal form, Young's orthogonal form, or Young's natural form. It would be interesting to understand how they are related to any of them.
\end{rmk*}

\begin{rmk*}
Moreover, the choice of these representations is not unique. Some of the freedom degrees that appear are not trivial ones such as reordering or negating coordinates. These are an effect of repeating eigenvalues, which account for most of the complications in the process of finding these matrices. The main challenge that arises from the preliminary results in this section, and is conjectured to be achievable, is as follows. 
\end{rmk*}

\begin{quest}
For every $k \in \mathbb{N}$, give a natural and explicit construction of representation matrices $R^{\lambda}(\sigma) \in {\mathbb{Q}}_k^{d_{\lambda} \times d_{\lambda}}$ for all $\sigma \in S_k$ and $\lambda \vdash k$, such that the matrix elements $\{\mathbf{R}_{ij}^{\lambda} \}$ fully diagonalize the asymptotic second-moment matrix of the normalized $k$-profile, as in Proposition~\ref{diagonal}.
\end{quest}

We mention a recent work that may be relevant to this direction. It is the spectral analysis of another problem on random permutations, the random-to-random shuffling, by Dieker and Saliola~\cite{dieker2018spectral}. 

\begin{rmk*} Another direction left for future research is studying the asymptotic behavior of the normalized profile beyond the second moment. This includes describing the non-normal limit distributions of the U-statistics that show up, and the joint distributions of several or all of them. 
\end{rmk*}

\begin{rmk*}
We include some technical information about the computer-assisted verification of Proposition~\ref{diagonal}. The verification takes 36 hours with SageMath version 8.1 on an AMD 2Ghz machine. Only the second moments of $\mathbf{P}_{6,11}$ and $\mathbf{P}_{6,12}$ were pre-prepared in another 36 hours, while distributed on up to 132 Intel 2.3Ghz processors running SageMath version 7.4.
\end{rmk*}

\section{Applications}
\label{apps}

Our viewpoint of the structure of the profile's distribution provides a unified framework for the null distribution of several nonparametric statistics, and throws light on other questions on random permutations. We describe some of these connections in this section.

\bigskip \noindent
\textbf{Rank Correlation Tests}
\nopagebreak \medskip \\ \noindent
Consider a sequence of $n$ independent paired samples $X_i = (Y_i, Z_i)$, drawn from a common continuous distribution $X_i \sim X$ on $\mathbb{R}^2$. Various real measures were suggested to detect, describe, and estimate the correlation between $Y$ and $Z$ based on these samples. Some of them are \emph{nonparametric}, using only the ranks of $\{Y_i\}$ and of $\{Z_i\}$. Since such a statistic is invariant under order-preserving maps for each coordinate, and under reordering the sequence of paired samples, it only depends on the induced permutation $\pi =\text{perm}(X_1,\dots,X_n)$, as defined in Section~\ref{back}. 

Kendall's tau coefficient~\cite{kendall1938new} is one such approach to quantify the correlation. In terms of permutations and patterns, let $\tau:S_n \to [-1,1]$,
$$ \tau(\pi) \;=\; \left\langle \mathbf{R}^{11}_{1,1} \;,\; \mathbf{P}_{2}(\pi) \right\rangle \;=\; P_{12}(\pi) - P_{21}(\pi) $$ 
Thus $\tau$ is a projection of the 2-profile along a matrix element, the one from the representation corresponding to $\lambda = (1,1)$.

Spearman's rho coefficient~\cite{spearman1904proof} is another nonparametric measure of correlation. Its standard definition uses the Pearson correlation between the ranks of the two coordinates of the samples. However, Spearman's rho is equivalently described using two matrix elements of the 3-profile:
\begin{align*}
\rho(\pi) \;=\; & \left\langle \tfrac{4n}{3n+3} \mathbf{R}^{21}_{2,2} - \tfrac{n-3}{3n+3} \mathbf{R}^{111}_{1,1} \;,\; \mathbf{P}_{3}(\pi) \right\rangle \\ \;=\; & P_{123} + P_{132} + P_{213} - P_{231} - P_{312} - P_{321} + O\left(\tfrac1n\right) 
\end{align*}
Here and below, our choice of matrix representations $R^{\lambda}$ is as in the previous section. The matrix element $\mathbf{R}^{21}_{2,2}$ corresponds to the first principal component of the limit distribution of $\sqrt{n}\,(\mathbf{P}_{3n} - \mathbf{u}_3)$.

The two coefficients $\tau$ and $\rho$ are often used in statistical hypothesis tests, to establish whether there is a significant correlation between $Y$ and $Z$. It is hence important to understand their distributions under the null hypothesis that $Y$ and $Z$ are independent, so that $\pi$ is uniform in~$S_n$. The null distributions of these tests are hence contained, as special cases, in our general discussion of the distribution of the profile of a random permutation.

The profile's decomposition provides tools to compare such nonparametric correlation tests. In terms of the 3-profile, we can express Kendall's $\tau$ using $\tfrac{8}{9} \mathbf{R}^{21}_{2,2} + \tfrac{1}{9} \mathbf{R}^{111}_{1,1}$ while for Spearman's $\rho$ we use $\tfrac{4}{3} \mathbf{R}^{21}_{2,2} - \tfrac{1}{3} \mathbf{R}^{111}_{1,1}$ asymptotically as $n \to \infty$. Assuming independence, both measures are asymptotically normal of order $1/\sqrt{n}$, being governed by the same term $\mathbf{R}^{21}_{2,2}$. The term~$\mathbf{R}^{111}_{1,1}$ makes a smaller, order-$1/n$ contribution. 

Instead of choosing between two tests, perhaps it makes sense to consider the whole family, 
$\alpha\, \mathbf{R}^{21}_{2,2} + (1-\alpha) \mathbf{R}^{111}_{1,1}$ for some range of real $\alpha$.
The choice of test may then depend on its power against various alternative models, of correlated data with typical noise and outliers.

\bigskip \noindent
\textbf{Circular Rank Correlation}
\nopagebreak \medskip \\ \noindent
Suppose that $X_i = (Y_i, Z_i)$ are independently sampled from a common continuous distribution on the torus~$S^1 \times S^1$. In this circular setting, it is natural to ask how much $Y$ and $Z$ tend to be related, approximately, by a homeomorphism between the two copies of~$S^1$. Such a relation can either preserve or reverse the orientation of the circle. 

Fisher and Lee~\cite{fisher1982nonparametric} suggested a nonparametric measure of correlation that depends on the \emph{circular ranks} of the coordinates of the samples. The ordering of three points in~$S^1$ can be either clockwise or counterclockwise. Let $\Delta:S_n\to[-1,1]$ be the proportion of triplets $\{i,j,k\}$ where $Y_i,Y_j,Y_k$ and $Z_i,Z_j,Z_k$ have the same ordering type, minus the proportion of triplets where they differ. This measure can be expressed in terms of permutations and patterns as
$$ \Delta(\pi) \;=\; \left\langle  \mathbf{R}^{111}_{1,1} \;,\; \mathbf{P}_{3}(\pi) \right\rangle \;=\; P_{123} - P_{132} - P_{213} + P_{231} + P_{312} - P_{321} $$
Here $\pi =\text{perm}(X_1,\dots,X_n)$ is defined by first removing one arbitrary point from each circle so that it maps to the real line. Note that the matrix element $\mathbf{R}^{111}_{1,1}$ belongs to the component~$V_2$ of $\mathbb{R}^{3!}$.

If one uses $\Delta$ as a test statistic, then under the null hypothesis that $Y$ and~$Z$ are independent, $\pi \in S_n$ is uniform, and $\Delta$ scales as $1/n$ by Theorem~\ref{proj1}. See~\cite{fisher1982nonparametric,  janson2015asymptotic, zeilberger2016doron} for more details on the non-normal limiting distribution of~$n\Delta$.

\begin{rmk*}
A crucial property of this projection of $\mathbf{P}_3$ is being symmetric under rotation of both the input or the output of the permutation. In general, restrictions to various invariant subspaces are useful in the study of the profile. Such symmetries include combinations of input rotation, output rotation, input reflection, output reflection, and inverting the permutation.
\end{rmk*}

\bigskip \noindent
\textbf{Rank Independence Tests}
\nopagebreak \medskip \\ \noindent
Let $X_1,\dots,X_n$ be independent paired samples from a common continuous distribution $X_i = (Y_i,Z_i) \sim X = (Y,Z)$ in $\mathbb{R}^2$ as above. Several statistics were suggested to detect \emph{any} kind of dependence between $Y$ and~$Z$. 

Hoeffding's independence test uses a nonparametric approach~\cite{hoeffding1948non}. Given the joint and marginal distribution functions $F_X$, $F_Y$, and~$F_Z$, let
$$ H \;=\; \iint \left(F_X(y,z) - F_Y(y)F_Z(z)\right)^2 dF_X(y,z) $$ 
For continuous distributions, $H=0$ if and only if $Y$ and $Z$ are independent. Blum, Kiefer and Rosenblatt~\cite{blum1961distribution} showed that replacing $dF_X(y,z)$ by $dF_Y(y)dF_Z(z)$ extends this condition to non-continuous distributions.

Hoeffding's test is based on a consistent estimator of~$H$ as a function of the $n$ samples. The functions $F_X$, $F_Y$, and $F_Z$ are replaced by the empirical distributions, which are cumulative counts of samples. Since $H$ contains products of up to five functions, this leads to a U-statistic with summation over five-sample subsets. In terms of the 5-profile, it has the following short description:
$$ D(\pi) \;=\; \tfrac{1}{30} \left\langle  \tfrac12 \mathbf{R}^{32}_{4,4} + \tfrac12 \mathbf{R}^{221}_{3,3} \;,\; \mathbf{P}_{5}(\pi) \right\rangle $$
A similar estimator for the BKR functional yields a U-statistic of order six, which also turns out to be expressible using the 5-profile:
$$ B(\pi) \;=\; \left\langle  \tfrac52 \mathbf{R}^{32}_{4,4} - \tfrac32 \mathbf{R}^{221}_{3,3} \;,\; \mathbf{P}_{5}(\pi) \right\rangle $$
Another variant of this test by Bergsma and Dassios~\cite{bergsma2010nonparametric, bergsma2014consistent} uses, up to constants and negligible $O(\tfrac1n)$ terms, the following sum of 4-pattern densities:
$$ P_{1234} + P_{1243} + P_{2134} + P_{2143} + P_{3412} + P_{3421} + P_{4312} + P_{4321} $$
This test is simply described in terms of the 4-profile or the 5-profile:
$$ {BD}(\pi) \;=\;  \left\langle \mathbf{R}^{22}_{2,2} \;,\; \mathbf{P}_{4}(\pi) \right\rangle \;=\; \left\langle \tfrac{9}{10} \mathbf{R}^{32}_{4,4} + \tfrac{1}{10} \mathbf{R}^{221}_{3,3} \;,\; \mathbf{P}_{5}(\pi) \right\rangle $$

Our analysis of the profile sheds new light on these independence tests. Their known $1/n$ scaling under the hypothesis of independence is a special case of Theorem~\ref{proj1}. The relative differences between them go to zero in probability, because the contribution of $\mathbf{R}^{221}_{3,3}$ has order~$n^{-3/2}$. Similar to the correlation tests above, it would be interesting to better understand the role of such lower-order terms in the tests' performance against various alternatives. The general form \mbox{$\alpha\,\mathbf{R}^{32}_{4,4} + (1-\alpha) \mathbf{R}^{221}_{3,3}$} contains  $B$, $D$, and ${BD}$ as special cases. 

\bigskip \noindent
\textbf{Quasirandom Permutations}
\nopagebreak \medskip \\ \noindent
Given a deterministic sequence $\pi_n \in S_n$ for every $n \in \mathbb{N}$, several equivalent notions of \emph{quasi-randomness} have been considered~\cite{cooper2004quasirandom}. These are defined as certain asymptotic properties of $\pi_n$, that hold for a random permutation with probability going to one. 

As noted in~\cite[App.~A]{hoeffding1948non}, the uniformity of $\text{perm}\left(X_1,X_2,X_3,X_4,X_5\right)$ in $S_5$ already implies independence, hence $X$ produces uniformly random permutations of any size. The same holds for $\text{perm}\left(X_1,X_2,X_3,X_4\right)$ by the consistency of the Bergsma--Dassios test. It follows by the arguments in \cite[cf.~\cite{kral2013quasirandom}]{hoppen2011testing}, that $\{\pi_n\}$ is quasirandom if and only if
$$\left\langle \mathbf{R}^{22}_{2,2} \;,\; \mathbf{P}_{4}(\pi_n) \right\rangle \;\to\;0 \;,$$
which hasn't been noted in those works. It would be interesting to characterize what other combinations of $k$-densities imply quasirandomness, and the profile's decomposition may play a role in that effort.

\bigskip \noindent
\textbf{Statistical Analysis of Rankings}
\nopagebreak \medskip \\ \noindent
Our treatment of the $k$-profile parallels the framework of spectral analysis of statistical data defined on non-abelian groups, as introduced by \mbox{Diaconis} \cite[Section 8B]{diaconis1988group}. An important practical issue that arises there is the arbitrary choice of Fourier bases, which might depend on matters of interpretation and convenience. 

The bases described here propose an answer for the case that the data takes values in~$S_k$. The matrix elements as in Section~\ref{diag} may give useful Fourier descriptions, especially in a setting where the samples may relate to occurrences of ordering types, or are possibly induced from a larger sequence, perhaps one with random-like features.

\bigskip \noindent
\textbf{Acknowledgments}
\nopagebreak \medskip \\ \noindent
The author thanks Eric Babson, Svante Janson, Nati Linial, and Doron Zeilberger for insightful discussions, and the anonymous
referees for their helpful comments.

\nopagebreak
Some of the computational work was carried out with the facilities of the  School of Computer Science and Engineering at HUJI, supported by
ERC~339096.

\bibliographystyle{alpha}
\bibliography{patterns}

\newcommand{\etalchar}[1]{$^{#1}$}
\begin{thebibliography}{KKRW15}

\bibitem[AAH{\etalchar{+}}02]{albert2002packing}
Michael~H Albert, Mike~D Atkinson, Chris~C Handley, Derek~A Holton, and Walter
  Stromquist.
\newblock On packing densities of permutations.
\newblock {\em Electron. J. Combin}, 9(1), 2002.

\bibitem[BD14]{bergsma2014consistent}
Wicher Bergsma and Angelos Dassios.
\newblock A consistent test of independence based on a sign covariance related
  to {K}endall's tau.
\newblock {\em Bernoulli}, 20(2):1006--1028, 2014.

\bibitem[Ber10]{bergsma2010nonparametric}
Wicher Bergsma.
\newblock Nonparametric testing of conditional independence by means of the
  partial copula.
\newblock 2010.

\bibitem[BH10]{burstein2010packing}
Alexander Burstein and Peter H{\"a}st{\"o}.
\newblock Packing sets of patterns.
\newblock {\em European Journal of Combinatorics}, 31(1):241--253, 2010.

\bibitem[BHL{\etalchar{+}}15]{balogh2015minimum}
J{\'o}zsef Balogh, Ping Hu, Bernard Lidick{\`y}, Oleg Pikhurko, Bal{\'a}zs
  Udvari, and Jan Volec.
\newblock Minimum number of monotone subsequences of length 4 in permutations.
\newblock {\em Combinatorics, Probability and Computing}, 24(04):658--679,
  2015.

\bibitem[BKR61]{blum1961distribution}
Julius~Robin Blum, Jack~Carl Kiefer, and Murray Rosenblatt.
\newblock Distribution free tests of independence based on the sample
  distribution function.
\newblock {\em The Annals of Mathematical Statistics}, pages 485--498, 1961.

\bibitem[B{\'o}n07]{bona2007copies}
Mikl{\'o}s B{\'o}na.
\newblock The copies of any permutation pattern are asymptotically normal.
\newblock {\em arXiv:0712.2792}, 2007.

\bibitem[B{\'o}n10]{bona2010three}
Mikl{\'o}s B{\'o}na.
\newblock On three different notions of monotone subsequences.
\newblock {\em Permutation Patterns}, 376:89--114, 2010.

\bibitem[B{\'o}n12]{bona2012combinatorics}
Mikl{\'o}s B{\'o}na.
\newblock {\em Combinatorics of Permutations}.
\newblock CRC Press, 2012.

\bibitem[Coo04]{cooper2004quasirandom}
Joshua~N Cooper.
\newblock Quasirandom permutations.
\newblock {\em Journal of Combinatorial Theory, Series A}, 106(1):123--143,
  2004.

\bibitem[Coo06]{cooper2006permutation}
Joshua~N Cooper.
\newblock A permutation regularity lemma.
\newblock {\em the Electronic Journal of Combinatorics}, 13(1):22, 2006.

\bibitem[Dia88]{diaconis1988group}
Persi Diaconis.
\newblock Group representations in probability and statistics.
\newblock {\em Lecture Notes -- Monograph Series}, 11:1--192, 1988.

\bibitem[DS18]{dieker2018spectral}
AB~Dieker and FV~Saliola.
\newblock Spectral analysis of random-to-random {M}arkov chains.
\newblock {\em Advances in Mathematics}, 323:427--485, 2018.

\bibitem[Eve18]{even2018gh}
Chaim Even{-Zohar}.
\newblock {\em GitHub repository \emph{patterns}}, 2018.
\newblock \url{http://github.com/chaim-e/patterns}.

\bibitem[FH91]{fulton1991representation}
William Fulton and Joe Harris.
\newblock {\em Representation Theory: A First Course}, volume 129.
\newblock Springer Science \& Business Media, 1991.

\bibitem[FL82]{fisher1982nonparametric}
NI~Fisher and AJ~Lee.
\newblock Nonparametric measures of angular-angular association.
\newblock {\em Biometrika}, pages 315--321, 1982.

\bibitem[Ful04]{fulman2004stein}
Jason Fulman.
\newblock Stein's method and non-reversible {M}arkov chains.
\newblock In {\em Stein's Method}, pages 66--74. Institute of Mathematical
  Statistics, 2004.

\bibitem[GGKK15]{glebov2015finitely}
Roman Glebov, Andrzej Grzesik, Tereza Klimo{\v{s}}ov{\'a}, and Daniel
  Kr{\'a}l'.
\newblock Finitely forcible graphons and permutons.
\newblock {\em Journal of Combinatorial Theory, Series B}, 110:112--135, 2015.

\bibitem[Has02]{hasto2002packing}
Peter~A Hasto.
\newblock The packing density of other layered permutations.
\newblock {\em Journal of Combinatorics}, 9(2):1, 2002.

\bibitem[HKM{\etalchar{+}}13]{hoppen2013limits}
Carlos Hoppen, Yoshiharu Kohayakawa, Carlos~Gustavo Moreira, Bal{\'a}zs
  R{\'a}th, and Rudini~Menezes Sampaio.
\newblock Limits of permutation sequences.
\newblock {\em Journal of Combinatorial Theory, Series B}, 103(1):93--113,
  2013.

\bibitem[HKMS11]{hoppen2011testing}
Carlos Hoppen, Yoshiharu Kohayakawa, Carlos~Gustavo Moreira, and Rudini~Menezes
  Sampaio.
\newblock Testing permutation properties through subpermutations.
\newblock {\em Theoretical Computer Science}, 412(29):3555--3567, 2011.

\bibitem[Hoe48a]{hoeffding1948class}
Wassily Hoeffding.
\newblock A class of statistics with asymptotically normal distribution.
\newblock {\em The Annals of Mathematical Statistics}, 19(3):293--325, 1948.

\bibitem[Hoe48b]{hoeffding1948non}
Wassily Hoeffding.
\newblock A non-parametric test of independence.
\newblock {\em The Annals of Mathematical Statistics}, pages 546--557, 1948.

\bibitem[Hof17]{hofer2017central}
Lisa Hofer.
\newblock A central limit theorem for vincular permutation patterns.
\newblock {\em arXiv preprint 1704.00650}, 2017.

\bibitem[Jan97]{janson1997gaussian}
Svante Janson.
\newblock {\em Gaussian Hilbert Spaces}, volume 129.
\newblock Cambridge University Press, 1997.

\bibitem[JNZ15]{janson2015asymptotic}
Svante Janson, Brian Nakamura, and Doron Zeilberger.
\newblock On the asymptotic statistics of the number of occurrences of multiple
  permutation patterns.
\newblock {\em Journal of Combinatorics}, 6:117--143, 2015.

\bibitem[KB13]{korolyuk2013theory}
Vladimir~S Korolyuk and Yu~V Borovskich.
\newblock {\em Theory of U-statistics}, volume 273.
\newblock Springer Science \& Business Media, 2013.

\bibitem[Ken38]{kendall1938new}
Maurice~G Kendall.
\newblock A new measure of rank correlation.
\newblock {\em Biometrika}, 30(1):81--93, 1938.

\bibitem[Kit11]{kitaev2011patterns}
Sergey Kitaev.
\newblock {\em Patterns in Permutations and Words}.
\newblock Springer Science \& Business Media, 2011.

\bibitem[KK14]{klimovsova2014hereditary}
Tereza Klimo{\v{s}}ov{\'a} and Daniel Kr{\'a}l'.
\newblock Hereditary properties of permutations are strongly testable.
\newblock In {\em Proceedings of the Twenty-Fifth Annual ACM-SIAM Symposium on
  Discrete Algorithms}, pages 1164--1173. Society for Industrial and Applied
  Mathematics, 2014.

\bibitem[KKRW15]{kenyon2015permutations}
Richard Kenyon, Daniel Kr{\'a}l', Charles Radin, and Peter Winkler.
\newblock Permutations with fixed pattern densities.
\newblock {\em arXiv:1506.02340}, 2015.

\bibitem[KP13]{kral2013quasirandom}
Daniel Kr{\'a}l' and Oleg Pikhurko.
\newblock Quasirandom permutations are characterized by 4-point densities.
\newblock {\em Geometric and Functional Analysis}, 23(2):570--579, 2013.

\bibitem[Lee90]{lee1990u}
Justin Lee.
\newblock {\em U-Statistics: Theory and Practice}, volume 110.
\newblock Marcel Dekker, Inc., 1990.

\bibitem[MT04]{marcus2004excluded}
Adam Marcus and G{\'a}bor Tardos.
\newblock Excluded permutation matrices and the {S}tanley--{W}ilf conjecture.
\newblock {\em Journal of Combinatorial Theory, Series A}, 107(1):153--160,
  2004.

\bibitem[Pri97]{price1997packing}
Alkes~Long Price.
\newblock Packing densities of layered *patterns, 1997.
\newblock Dissertations available from ProQuest. AAI9727276.

\bibitem[PS10]{presutti2010packing}
Cathleen~Battiste Presutti and Walter Stromquist.
\newblock Packing rates of measures and a conjecture for the packing density of
  2413.
\newblock {\em Permutation Patterns}, 376:287--316, 2010.

\bibitem[{Sag}18]{sagemath}
{Sage Developers}.
\newblock {\em {S}ageMath, the {S}age {M}athematics {S}oftware {S}ystem,
  {V}ersions 7.4 and 8.1}, 2018.
\newblock \url{http://www.sagemath.org}.

\bibitem[Spe04]{spearman1904proof}
Charles Spearman.
\newblock The proof and measurement of association between two things.
\newblock {\em The American Journal of Psychology}, 15(1):72--101, 1904.

\bibitem[SS17]{sliacan2017improving}
Jakub Sliacan and Walter Stromquist.
\newblock Improving bounds on packing densities of 4-point permutations.
\newblock {\em arXiv:1704.02959}, 2017.

\bibitem[Zei16]{zeilberger2016doron}
Doron Zeilberger.
\newblock Doron {G}epner's statistics on words in $\{$1, 2, 3$\}$ is (most
  probably) asymptotically logistic.
\newblock {\em arXiv:1604.00663}, 2016.

\end{thebibliography}

\end{document}